\documentclass[10pt, reqno]{amsart}

\usepackage{graphicx}
\usepackage{times}

\usepackage{xcolor}
\usepackage{hyperref}
\hypersetup{
    linktoc=all,
    colorlinks=true
}
\usepackage{stmaryrd}
\hypersetup{
    colorlinks = true,
    linkcolor = blue,
    anchorcolor = blue,
    citecolor = blue,
    filecolor = blue,
    urlcolor = blue
    }

\newtheorem{theorem}{Theorem}[section]

\newtheorem{lemma}[theorem]{Lemma}
\newtheorem{corollary}[theorem]{Corollary}

\newtheorem{notation}[theorem]{Notation}

\theoremstyle{definition}
\newtheorem{definition}[theorem]{Definition}

\theoremstyle{remark}
\newtheorem{remark}[theorem]{Remark}
\numberwithin{equation}{section}

\usepackage[top=0.5in, bottom=0.75in, left=0.75in, right=0.75in]{geometry}

\usepackage[scr]{rsfso}
\usepackage[english]{babel}
\usepackage{fancyhdr}
\usepackage{amsmath}
\usepackage{amsthm}
\usepackage{amssymb}
\usepackage{eucal}
\usepackage{comment}
\usepackage{enumitem}
\setlist{leftmargin=*}
\usepackage[integrals]{wasysym}
\newcommand\nc{\newcommand}
\nc{\on}{\operatorname}
\nc{\E}{\mathbf{E}}
\nc{\R}{\mathbb R}
\nc{\C}{\mathbb C}
\nc{\Q}{\mathbb Q}
\nc{\Z}{\mathbb Z}
\nc{\N}{\mathbb N}
\nc{\F}{\mathbb F}
\nc{\wt}{\widetilde}
\nc{\ol}{\overline}
\nc{\short}[3]{0 \longrightarrow #1 \longrightarrow #2 \longrightarrow #3 \longrightarrow 0}
\nc{\pd}[2]{\frac{\partial #1}{\partial #2}}
\nc{\rnc}{\renewcommand}
\nc{\e}{\varepsilon}
\nc{\DMO}{\DeclareMathOperator}
\nc{\grad}{\nabla}
\nc{\Exp}{\mathbf{Exp}}
\nc{\fsp}{\fontdimen2\font=2.0pt}

\rnc{\leq}{\leqslant}
\rnc{\geq}{\geqslant}
\rnc{\d}{\mathrm{d}}
\rnc{\O}{\mathrm{O}}
\rnc{\exp}{\mathbf{Exp}}
\newenvironment{nouppercase}{%
  \renewcommand{\uppercasenonmath}[1]{}}{}
\title{\Large Stochastic Burgers Equation via Energy Solutions \\ from Non-Stationary Particle Systems}
\author{\large Kevin Yang}
\usepackage{setspace}
\begin{document}
\setstretch{1.0}
\fontdimen2\font=2.17pt
\raggedbottom
\begin{nouppercase}
\maketitle
\end{nouppercase}

\begin{abstract}
\fsp We prove that the stochastic Burgers equation, which is related to the Kardar-Parisi-Zhang/KPZ equation via weak derivative, is a ``critical" scaling limit for density fluctuations for a family of non-integrable and non-stationary interacting particle systems. The models we consider cannot be linearized by a microscopic Cole-Hopf transform or studied directly by the existing energy solution theory of \cite{GJ15}. We develop a novel method based on comparison to stationary models, a technique that has not yet been applied to universality of the KPZ equation and nontrivially expands the set of models for which universality is confirmed. We also study \emph{crossover fluctuations} and prove, for the first time, the full transition and phase diagram from Gaussian to KPZ fluctuations for non-stationary interacting particle systems, which has not even been done yet for integrable models. We emphasize the method developed herein applies to a general class of models/particle systems, but we restrict to a class of zero-range systems whose non-stationary versions have received widespread interest but had not been treated in the context of KPZ until now as well as class of non-simple exclusion processes that we comment on.
\end{abstract}

\section{Introduction}
This paper addresses the universality problem of the \emph{Kardar-Parisi-Zhang}/\emph{KPZ} stochastic PDE. This equation was derived in \cite{KPZ} to provide a \emph{universal} description of interface fluctuation statistics, though rigorous proof has been limited to a handful of models that enjoy rich algebraic/solvability structure or models that enjoy both a sufficiently explicit invariant measure and starting at or extremely close to said invariant measure. We additionally consider the \emph{stochastic Burgers equation}/\emph{SBE}, whose solution is related to that of KPZ via $\mathbf{u}=\grad\mathbf{h}$. Formally, these equations are written below in which $\alpha,\zeta$ are positive and $\beta\in\R$:
\begin{align}
\partial_{t}\mathbf{h} \ = \ 2^{-1}\alpha\Delta\mathbf{h} + \beta(\grad\mathbf{h})^{2} + \zeta\xi \quad \mathrm{and} \quad \partial_{t}\mathbf{u} \ = \ 2^{-1}\alpha\Delta\mathbf{u} + \beta\grad(\mathbf{u}^{2}) + \zeta\grad\xi. \label{eq:SPDE}
\end{align}
The SPDEs \eqref{eq:SPDE} are both posed on the space-time $(t,x)\in\R_{\geq0}\times\mathbb{T}$, in which $\mathbb{T}=\R/\Z$ is the unit-length one-dimensional torus. For convenience, we refer to the SPDEs in \eqref{eq:SPDE} as $\mathrm{KPZ}(\alpha,\beta,\zeta)$ and $\mathrm{SBE}(\alpha,\beta,\zeta)$, respectively. In the case $\beta=0$, we adopt the notation $\mathrm{EW}(\alpha,\zeta)=\mathrm{KPZ}(\alpha,0,\zeta)$ and $\grad\mathrm{EW}(\alpha,\zeta)=\mathrm{SBE}(\alpha,0,\zeta)$, which denote the solution to SPDEs that are often referred to as the \emph{Edwards-Wilkinson limit} or \emph{Additive Stochastic Heat Equation} or \emph{ASHE} and its gradient that itself is referred to as a \emph{conservative ASHE}. When we refer to these equations without specific $\alpha,\beta,\zeta$, we write $\mathrm{KPZ}$ or $\mathrm{SBE}$ or $\mathrm{EW}$ or $\grad\mathrm{EW}$.

Progress towards the aforementioned \emph{universality} of SPDEs in \eqref{eq:SPDE} has been fairly modest, especially in the context of density fluctuations, or any other observables, of interacting particle systems. In the seminal work \cite{BG}, the authors proved that the \emph{integrated} density fluctuation, which corresponds to the corner growth model of height functions, of the asymmetric simple exclusion process/ASEP converges to $\mathrm{KPZ}$ with explicit $\alpha,\beta,\zeta$ parameters. Several papers including \cite{CST, CT} followed work of \cite{BG} and established $\mathrm{KPZ}$ scaling limits. However, these rely heavily on an incredibly rare and special \emph{integrability} property of the models addressed therein. On the other hand, in \cite{GJ15} an alternative \emph{energy solution theory} for $\mathrm{SBE}$ was developed. Together with a key identification in \cite{GP} of energy solution theory with other notions of solutions of KPZ and SBE, the energy solution method of \cite{GJ15} proved $\mathrm{SBE}$ scaling limits for density fluctuations in a wide class of interacting particle systems; see \cite{GJS15}. However, the energy solution theory has only been successfully applied to models with an explicit invariant measure; these models are also assumed to start at such an invariant measure, or with bounded relative entropy with respect to one. Provided many of the interesting connections between KPZ and random matrix theory, for example, happen well-beyond invariant measures, this motivates proving the universality of $\mathrm{KPZ}$ and $\mathrm{SBE}$ for non-stationary models, especially non-integrable models. We make progress in this direction by applying the energy solution theory for a nontrivial family of \emph{non-stationary} models that are by no means perturbations of the integrable models in the aforementioned papers \cite{BG, CST, CT}. This provides the long sought-after first application of energy solution theory for non-stationary systems. Moreover, our application of energy solution theory confirms the complete phase diagram for scaling limits of fluctuation fields that was shown for stationary models in \cite{GJ15, GJS15}, showing the sharp transition from $\mathrm{EW}$/$\grad\mathrm{EW}$ scaling limits to $\mathrm{KPZ}$/$\mathrm{SBE}$ scaling limits. This has not been accomplished for non-stationary models prior to this work, even for integrable models in aforementioned papers \cite{BG, CST, CT}. As for the models for which our method holds, we have in mind a fairly large set of particle systems. However, for clearer exposition, we specialize to zero-range models and comment at the end of this section on a general class of non-simple exclusion processes.

Let us expand on the universality/scaling limit problem discussed above in more detail. Ultimately, the two approaches via integrability and via energy solution theory are each based on making sense of the singular nonlinearity in $\mathrm{KPZ}$ and $\mathrm{SBE}$ in a way that can be realized at the level of particle systems. The two approaches are explained below.
\begin{itemize}
\item In \cite{BG}, Bertini and Giacomin introduce the following (multiplicative noise) \emph{stochastic heat equation}/\emph{SHE}. Its solution $\mathbf{z}$ \emph{defines} a solution to $\mathrm{KPZ}$ via $\mathbf{h}=\lambda^{-1}\log\mathbf{z}$, where $\lambda=\beta\alpha^{-1}$. In particular, the following SHE on $\R_{\geq0}\times\mathbb{T}$ admits solutions that are continuous in space-time and strictly positive with probability 1 if the initial data is non-zero and non-negative. This is the \emph{comparison principle} of \cite{Mu}. We now write the SHE; we refer to its solution as $\mathrm{SHE}(\alpha,\beta,\zeta)$ or $\mathrm{SHE}$ more generally:
\begin{align}
\partial_{t}\mathbf{z} \ = \ 2^{-1}\alpha\Delta\mathbf{z} + \beta\alpha^{-1}\zeta\mathbf{z}\xi. \label{eq:SHE}
\end{align}
A lot of the integrability of $\mathrm{KPZ}$ is clarified via the SHE; see \cite{C11}. Moreover, the models in \cite{BG, CST, CT} all satisfy a special property. Exponentiating their respective height functions provides an exact discrete analog of \eqref{eq:SHE} without error terms, so the authors may essentially follow their noise to uncovering $\mathrm{KPZ}$ scaling limits. We clarify the transform $\mathbf{z}=\exp(\lambda\mathbf{h})$ is often referred to as the \emph{Cole-Hopf transform}, and it implicitly handles the singular behavior of $\mathrm{KPZ}$, which we explain more of in the following bullet point, by exploiting an important \emph{integrability} of the KPZ equation.
\item In \cite{GJ15}, Goncalves and Jara introduce the energy solution theory for $\mathrm{SBE}$ that is basically a nonlinear martingale problem formulation for $\mathrm{SBE}$. Provided any test function $\Phi$ on the unit torus, the one-dimensional stochastic process defined below is required to be a martingale with the correct quadratic variation $\zeta^{2}\|\grad\Phi\|_{\mathscr{L}^{2}(\mathbb{T})}^{2}t$, thus equal to the pairing of $\grad\xi$ with $\Phi$ in law:
\begin{align}
\int_{\mathbb{T}}\mathbf{u}(t,x)\Phi(x)\d x - \int_{\mathbb{T}}\mathbf{u}(0,x)\Phi(x)\d x - \int_{0}^{t}\int_{\mathbb{T}}\mathbf{u}(s,x)\cdot \left(\alpha\Delta\Phi(x)\right)\d x\d s + \int_{0}^{t}\int_{\mathbb{T}}|\mathbf{u}(s,x)|^{2}\cdot \left(\beta\grad\Phi(x)\right)\d x\d s.
\end{align}
The energy solution theory of \cite{GJ15} is thus a weak formulation, in both an analytic and probabilistic sense, of $\mathrm{SBE}$. However, because $\xi$ is a rough noise term, the solution $\mathbf{u}$ is not an actual function but instead a generalized function, so its square is not classically or canonically defined. Remedying this singular behavior is the technical goal of energy solution theory. For this purpose it requires strong assumptions on the Gaussian invariant measure/initial measure of $\mathrm{SBE}$ when viewed as a Markov process on the space of generalized functions. As a lot of the previous literature on scaling limits for density fields/fluctuations of interacting particle systems has been built on martingale problems and weak formulations of PDEs and SPDEs, in \cite{GJ15,GJS15} the authors are able to establish $\mathrm{SBE}$ scaling limits for density fluctuations that exhibit the necessary invariant measure details needed to make sense of the nonlinearity in $\mathrm{SBE}$. They are also able to confirm the aforementioned phase diagram as well.
\end{itemize}
As for deriving KPZ and SBE limits for interacting particle systems in this paper, the general philosophy is similar to the two bullet points above in that we look for features of the models considered herein shared by the proposed scaling limit SPDEs. However, the method itself is novel. The only additional ingredient that we need is a \emph{monotonicity}-type property for the height functions of these models, which we realize as a microscopic version of the important comparison principle for $\mathrm{SHE}$ of \cite{Mu}. One broad interpretation of our work is that this microscopic comparison principle is one method to deal with the singular features of $\mathrm{SBE}$ that were dealt with via integrability or equilibrium calculations in \cite{BG, CST, CT, GJ15, GJS15}. The main technical challenge in this work is to \emph{approximate} general non-stationary models by initially ``almost-stationary" ones, namely interpolation between stationary models and general non-stationary models that lets us compare the latter to those for which the energy solution theory of \cite{GJ15, GJS15} is applicable. In particular, comparison with stationary data, for which the quadratic KPZ nonlinearity is well-defined, will replace the necessity in addressing the singular nonlinearity in $\mathrm{SBE}$ in this paper. Such comparison method has not been introduced or touched on in the context of universality of $\mathrm{KPZ}$/$\mathrm{SBE}$. In this paper we will use it to take a significant step towards universality and prove $\mathrm{KPZ}$/$\mathrm{SBE}$ scaling limits for density fluctuations of an interacting particle system with possibly unbounded number of particles per point that is not just a perturbation of ASEP. We hope that our method, provided both its youth and the fact that microscopic comparison principles, and their variations, are shared among a much larger set of interacting particle systems than integrability, gives another avenue towards attacking the long-standing problem of universality for $\mathrm{KPZ}$/$\mathrm{SBE}$.

To fix ideas, we fix a specific interacting particle system to analyze in this paper; we introduce a set of non-simple exclusion processes at the end of this introduction to comment on general applicability of our methods. The following is an \emph{attractive} zero-range process with a technical condition to prevent blow-up and an important spectral gap. We provide the example $\mathfrak{g}(k)=k$ to clarify the construction in Definition \ref{definition:intro1} below.
\begin{definition}\label{definition:intro1}
\fsp Define the discrete torus $\mathbb{T}_{N}=\Z/N\Z$. Provided any subset $\mathbb{I}\subseteq\mathbb{T}_{N}$, we define $\Omega_{\mathbb{I}}=\Z_{\geq0}^{\mathbb{I}}$ as a configuration space of particles on $\mathbb{I}$. For convenience, we also adopt the convention $\Omega=\Omega_{\mathbb{T}_{N}}$. Elements in $\Omega$ will be denoted by $\eta$, and for any $x\in\mathbb{T}_{N}$ we let $\eta_{x}\in\Z_{\geq0}$ indicate the number of particles at $x\in\mathbb{T}_{N}$. We now consider a function $\mathfrak{g}:\Omega\to\R_{\geq0}$ with non-negative values that satisfies the following list of properties in which $\sigma$ is positive and \emph{fixed} throughout this paper:
\begin{itemize}
\item We have $\mathfrak{g}(0) = 0$ and the Lipschitz property $\sup_{k\geq0}|\mathfrak{g}(k+1)-\mathfrak{g}(k)|<\infty$. From this, we construct a one-parameter class $\{\mathbf{P}_{\alpha}\}_{\alpha>0}$ of product probability measures on $\Omega$ whose one-dimensional marginals are all defined by the following Poisson-type formula, in which $\mathbf{Z}_{\alpha}$ is a normalization:
\begin{align}
\mathbf{P}_{\alpha}(k) \ = \ \mathbf{Z}_{\alpha}^{-1}\alpha^{k}{\prod}_{j=1}^{k}\mathfrak{g}(j)^{-1}.
\end{align}
\item Define $\alpha_{\sigma}$ to be the unique positive constant so that the $\sigma$-density is achieved, namely that $\E_{\alpha_{\sigma}}\eta_{0}=\sigma$. In turns out $\alpha_{\sigma}$ is unique, if it exists, because of monotonicity in the previous expectation with respect to $\alpha$; see Section 2.3 of \cite{GJS15}. Here, the $\E_{\alpha}$ expectation is with respect to $\mathbf{P}_{\alpha}$ and $\eta_{0}\sim\E_{\alpha}$. We now state the assumption of this bullet point. We assume there exists $\delta$ positive so that $\E_{\alpha_{\sigma}}\eta_{0}^{2+\delta}$ is finite.
\item The function $\mathfrak{g}$ is monotone non-decreasing, so that for any $k \in \Z_{\geq 0}$, we have $\mathfrak{g}(k+1) \geq \mathfrak{g}(k)$.
\item There exists $\e_{0}\in\R_{>0}$ and $x_{0}\in\Z_{\geq0}$ such that $\inf_{k\in\Z_{\geq0}} (\mathfrak{g}(k+x_{0})-\mathfrak{g}(k)) \geq \e_{0}$. This implies an $\e_{0},x_{0}$-dependent spectral gap with respect to $\mathbf{P}_{\alpha}$ for the $\mathscr{L}_{N}$ generator to be defined below, whose symmetric (antisymmetric) part is $\mathscr{L}_{N,\mathrm{S}}$ ($\mathscr{L}_{N,\mathrm{A}}$).
\end{itemize}
We now specify the infinitesimal generator of the zero-range process of interest in this paper. First, we let $\mathscr{L}_{x,\pm}$ denote the speed-1 infinitesimal generator for the process that moves a particle $x\to x\pm1$, where $x\in\mathbb{T}_{N}$ is arbitrary. Provided an asymmetry speed $\mathfrak{a}\in\R$, we define the infinitesimal generator $\mathscr{L}_{N}=\mathscr{L}_{N,\mathrm{S}}+\mathfrak{a}\mathscr{L}_{N,\mathrm{A}}$, where we take $\gamma\geq1/2$ below:
\begin{align}
\mathscr{L}_{N,\mathrm{S}} = 2^{-1}N^{2}{\sum}_{x\in\mathbb{T}_{N}}\mathfrak{g}(\eta_{x})\mathscr{L}_{x,+} + 2^{-1}N^{2}{\sum}_{x\in\mathbb{T}_{N}}\mathfrak{g}(\eta_{x+1})\mathscr{L}_{x+1,-} \quad \mathrm{and} \quad \mathscr{L}_{N,\mathrm{A}} = 2^{-1}N^{2-\gamma}{\sum}_{x\in\mathbb{T}_{N}}\mathfrak{g}(\eta_{x})\mathscr{L}_{x,+}. \nonumber
\end{align}
We let $\eta^{N}_{T}$ be the process at time $T$ under the $\mathscr{L}_{N}$ process and let $\eta_{T,x}^{N}\in\Z_{\geq0}$ be the number of particles at $x\in\mathbb{T}_{N}$ at time $T$.
\end{definition}
The particle system described by Definition \ref{definition:intro1} is an attractive zero-range process. It has a one-parameter family of probability measures $\mathbf{P}_{\alpha}$, which turn out to be invariant measures for the $\mathscr{L}_{N}$ generator, that exhibit a uniformly bounded, from below, spectral gap. It turns out this zero-range process satisfies a ``gradient condition" as was assumed and specified in Section 2.3 of \cite{GJS15}. These three properties of the process will basically be all we need, as long as we have additional technical conditions such as no blow-up and we have $2+\delta$-moments of the particle density under the aforementioned product measure $\mathbf{P}_{\alpha_{\sigma}}$.
\begin{remark}\label{remark:intro2}
\fsp If for our choice of $\sigma$ in Definition \ref{definition:intro1}, there exists no $\alpha_{\sigma}$ that is required in Definition \ref{definition:intro1}, we will ignore this value of $\sigma$. In particular, we will not consider such values of $\sigma$ in this paper, following Section 2.3 of \cite{GJS15} for zero-range processes therein.
\end{remark}
\subsection{Density Fluctuations and Height Function}
We spend the current subsection introducing the main objects of interest. The first of these is the density fluctuation field that will ultimately converge to $\mathrm{SBE}$ in the main result of this paper. In principle, the second height function converges to $\mathrm{KPZ}$ as well, though in this paper we will mostly use it as a technical integrated version of the density fluctuation field; see Lemma \ref{lemma:SBP}.
\begin{definition}\label{definition:intro3}
We let the \emph{density fluctuation field} with particle density parameter $\sigma \in \R_{>0}$ be the random space-time field $\mathbf{Y}_{T}^{N}: \mathscr{C}^{\infty}(\mathbb{T}) \to \R$ defined by the following action for any test function $\Phi \in \mathscr{C}^{\infty}(\mathbb{T})$:
\begin{align}
\mathbf{Y}_{T}^{N}(\Phi) \ \overset{\bullet}= \ N^{-\frac12}{\sum}_{X \in \mathbb{T}_{N}} \left( \eta_{T,X}^{N} - \sigma \right) \cdot \Phi(N^{-1}X - c_{N,\mathfrak{a},\gamma,\sigma}T).
\end{align}
We now define the following characteristic speed in $\Phi$ from the $\mathbf{Y}^{N}$ definition above. In what follows, we recall $\mathfrak{a}$ is the asymmetry introduced in Definition \ref{definition:intro1} for the particle system, and $\gamma$ controls the asymptotic strength of this asymmetry. We clarify $\eta_{0}$ on the RHS below is distributed via $\E_{\alpha_{\sigma}}$ introduced in Definition \ref{definition:intro1}. We are thereby differentiating the expectation on the RHS below with respect to its dependence on $\sigma$ and then evaluating this derivative at our fixed choice of $\sigma$ in Definition \ref{definition:intro1}:
\begin{align}
c_{N,\mathfrak{a},\gamma,\sigma} \ &\overset{\bullet}= \ 2^{-1}\mathfrak{a} N^{-\gamma+1}\partial_{\sigma}\left(\E_{\alpha_{\sigma}}\mathfrak{g}(\eta_{0})\right)|_{\sigma}.
\end{align}
\end{definition}
\begin{definition}\label{definition:HF}
First, provided any $T \in \R_{\geq0}$, define $\mathbf{h}_{T,0}^{N}$ as $N^{-\frac12}$ times the net flux of particles across the periodic boundary $N \rightsquigarrow 1$ in the time-interval $[0,T] \subseteq \R_{\geq 0}$, with the left-ward moving particles going $1 \rightsquigarrow N$ counting as $+1$ net flux and right-ward moving particles going $N \rightsquigarrow 1$ counting as $-1$ net flux. For $(T,X) \in \R_{\geq 0} \times \mathbb{T}_{N}$, we define the height function below in which $Y<X$ on the torus $\mathbb{T}_{N}$ means the inherited ordering $Y<X$ when identifying $\mathbb{T}_{N} \simeq [1,N]\cap\Z$.
\begin{align}
\mathbf{h}_{T,X}^{N} \ \overset{\bullet}= \ \mathbf{h}_{T,0}^{N} \ + \ N^{-\frac12}{\sum}_{Y < X} \left( \eta_{T,Y}^{N} - \sigma \right).
\end{align}
\end{definition}
\subsection{Main Results}
We now present the main theorems of this paper, which we recall are universality/scaling limit results for $\mathrm{KPZ}$/$\mathrm{SBE}$ as well as a phase diagram showing transition to these scaling limits from $\mathrm{EW}$/$\grad\mathrm{EW}$ scaling limits, all in the context of the zero-range processes of this paper. We first introduce the main type of initial conditions of this paper. We clarify the following does not refer to the $\mathbf{P}_{\alpha}$ invariant measures. We also clarify the following initial data assumes a \emph{deterministic} limit. We comment on how to remove this deterministic assumption towards the end of this section.
\begin{definition}\label{definition:IC}
\fsp We say a sequence of probability measures $\{\mathbf{P}_{N}\}_{N\in\Z_{>0}}$ on $\Omega_{\mathbb{T}_{N}}$ is \emph{deterministic and near-stationary} if under these measures, there is a deterministic continuous function $\mathbf{h}^{\infty}$ such that $\mathbf{h}_{0,N\bullet}^{N} \to_{N \to \infty} \mathbf{h}^{\infty}_{\bullet}$ uniformly on $[0,1]$ in probability.
\end{definition}
One example of deterministic and near-stationary initial measures are ones that are supported on the configuration where $\eta_{x}$ is equal to $x$ mod 2. In this case, taking $\sigma=2^{-1}$, we have the limit $\mathbf{h}^{N}_{0,N\bullet}\to0$. This is the \emph{flat} measure; see \cite{C11}. 

We now begin presenting our main results, the first of which is a central limit theorem for the dynamical fluctuations $\mathbf{Y}^{N}$ where the scaling limit is $\grad\mathrm{EW}$. Throughout this subsection and results below, we will refer to the Skorokhod topology on the Skorokhod space $\mathscr{D}(\R_{\geq0},\mathscr{D}(\mathbb{T}))$ of cadlag paths valued in the topological dual $\mathscr{D}(\mathbb{T})$ of the space $\mathscr{C}^{\infty}(\mathbb{T})$ of smooth functions on $\mathbb{T}$; see \cite{Bil}, for example. We emphasize that the constants appearing in the SPDE limits below can be found in (2.11) of \cite{GJS15} with details specific to zero-range processes in Section 2.3 of \cite{GJS15}. First, we record Theorem 2.2 in \cite{GJS15} for stationary interacting particle systems.
\begin{theorem}\label{theorem:ASHEZRPEq}
\fontdimen2\font=1.7pt Suppose that the interacting particle system is initially distributed via the grand-canonical invariant measure $\mathbf{P}_{\alpha_{\sigma}}$; this is defined in \emph{Definition \ref{definition:intro1}}. Also, assume $\gamma > \frac12$. The space-time process $\mathbf{Y}^{N}$ converges, in law, in the Skorokhod space $\mathscr{D}(\R_{\geq 0},\mathscr{D}(\mathbb{T}))$ to the solution of $\grad\mathrm{EW}(\alpha,\zeta)$ with the following parameters and the initial condition given by a spatial Gaussian white noise of mean zero and variance $\E_{\alpha_{\sigma}}\eta_{0}^{2}$, where $\E_{\alpha_{\sigma}}$ is expectation with respect to $\mathbf{P}_{\alpha_{\sigma}}$ and $\eta_{0}\sim\mathbf{P}_{\alpha_{\sigma}}$:
\begin{subequations}
\begin{align}
\alpha \ = \ \partial_{\sigma}\left(\E_{\alpha_{\sigma}}\mathfrak{g}(\eta_{0})\right)|_{\sigma} \quad \mathrm{and} \quad \zeta \ = \ \left(\E_{\alpha_{\sigma}}\mathfrak{g}(\eta_{0})\right)^{1/2}.
\end{align}
\end{subequations}
The same is true if instead of the grand-canonical invariant measure of density parameter $\sigma \in \R_{>0}$, the particle system is initially distributed via any sequence of measures $\{\mathbf{P}_{N}\}_{N\in\Z_{>0}}$ on $\Omega_{\mathbb{T}_{N}}$ which is uniformly bounded in relative entropy with respect to the aforementioned grand-canonical ensembles. The initial data of the limit $\grad\mathrm{EW}(\alpha,\zeta)$ is now the limit of the functional $\mathbf{Y}^{N}_{0}$ under $\mathbf{P}_{N}$ measures. The relative entropy is defined in \emph{\eqref{eq:RE}}.
\end{theorem}
The analog of Theorem \ref{theorem:ASHEZRPEq} for near-stationary initial conditions is the first main result, obtained by comparing to the stationary case in Theorem \ref{theorem:ASHEZRPEq}. We emphasize $\grad\mathrm{EW}$ has a direct and elementary solution theory for any continuous initial condition.
\begin{theorem}\label{theorem:ASHEZRP}
The results of \emph{Theorem \ref{theorem:ASHEZRPEq}} are true if the grand-canonical ensemble is replaced by any sequence of deterministic and near-stationary initial conditions on $\Omega_{\mathbb{T}_{N}}$. The initial data of the scaling limit $\grad\mathrm{EW}(\alpha,\zeta)$ is $\partial_{X}\lim_{N\to\infty}\mathbf{h}_{0,N\bullet}^{N}=\partial_{X}\mathbf{h}^{\infty}_{\bullet}$.
\end{theorem}
Theorem \ref{theorem:ASHEZRPEq} and Theorem \ref{theorem:ASHEZRP} are known as \emph{crossover fluctuations} as they describe the continuum model for $\mathbf{Y}^{N}$ for sufficiently weak asymmetry $\gamma \in (\frac12,\infty]$, and at the critical value $\gamma = \frac12$, we observe a phase transition away from the conservative ASHE and to the SBE. Again, we present the following result, which is Theorem 2.3 in \cite{GJS15}, which holds only with respect to the invariant measures and their bounded relative entropy perturbations as initial conditions/measures.
\begin{theorem}\label{theorem:SBEZRPEq}
\fontdimen2\font=1.7pt Suppose that the interacting particle system is initially distributed via the grand-canonical invariant measure $\mathbf{P}_{\alpha_{\sigma}}$ of density parameter $\sigma \in \R_{>0}$. Assume $\gamma = \frac12$. The results of \emph{Theorem \ref{theorem:ASHEZRPEq}} hold upon replacing $\grad\mathrm{EW}(\alpha,\zeta)$ with $\mathrm{SBE}(\alpha,\beta,\zeta)$ with parameters $\alpha,\zeta$ in \emph{Theorem \ref{theorem:ASHEZRPEq}} and
\begin{align}
\beta = 2^{-1}\mathfrak{a}\partial_{\sigma}^{2}\left(\E_{\alpha_{\sigma}}\mathfrak{g}(\eta_{0})\right)|_{\sigma}.
\end{align}
\end{theorem}
We refer to the discussion prior to Theorem 2.3 in \cite{GJS15} for a precise definition of solving $\mathrm{SBE}(\alpha,\beta,\zeta)$ using energy solution theory. The analog of Theorem \ref{theorem:SBEZRPEq} for deterministic and near-stationary initial conditions is the following result, again obtained by comparing with the stationary case. For the following, however, by solution to $\mathrm{SBE}(\alpha,\beta,\zeta)$, we mean the Cole-Hopf solution $\partial_{X}\lambda^{-1}\log\mathbf{z}$ defined prior to \eqref{eq:SHE}, which we emphasize agrees (see Theorem 2.13 in \cite{GP}) with the aforementioned energy solution from \cite{GJS15} for stationary models and for models which have uniformly bounded relative entropy with respect to stationary models \cite{GP17}. We emphasize here that energy solution theory of \cite{GJ15} provides no canonical notion of solution to $\mathrm{SBE}$ for general non-stationary initial data, whereas the Cole-Hopf solution is well-defined not just for Holder initial data but any possibly random continuous initial data (\cite{ACQ} obtains well-posedness for Dirac initial data, which easily extends via superposition of Dirac initial data to continuous initial data, since the $\mathrm{SHE}$ is linear).
\begin{theorem}\label{theorem:SBEZRP}
The results of \emph{Theorem \ref{theorem:SBEZRPEq}} are true if the grand-canonical ensemble is replaced by any sequence of deterministic and near-stationary initial conditions on $\Omega_{\mathbb{T}_{N}}$. The initial condition of the limit $\grad\mathrm{EW}(\alpha,\zeta)$ is $\partial_{X}\lim_{N\to\infty}\mathbf{h}_{0,N\bullet}^{N}=\partial_{X}\mathbf{h}^{\infty}_{\bullet}$.
\end{theorem}
\subsection{Additional Background and Comments}
Let us first comment on the \emph{deterministic} nature of the initial data considered in Theorem \ref{theorem:ASHEZRP} and Theorem \ref{theorem:SBEZRP}, in particular how to remove this assumption for another ``more common" assumption. In \cite{BG, CST, CT}, the height function is initially assumed to be basically Holder continuous of strictly positive Holder exponent. To deduce Theorem \ref{theorem:ASHEZRP} and Theorem \ref{theorem:SBEZRP} for these possibly non-deterministic initial conditions, because the unit torus $\mathbb{T}$ is compact, the embedding of Holder functions on $\mathbb{T}$, with positive exponent, into the space of continuous functions on $\mathbb{T}$ is compact. Thus, within any arbitrarily small but positive and fixed scale $r$, we approximate the potentially random limit of initial height functions by a \emph{finite} $r$-dependent number of deterministic and near-stationary height functions because of compactness. Our proofs of Theorems \ref{theorem:ASHEZRP} and \ref{theorem:SBEZRP}, which prove uniform-in-$N$ continuity of the density fluctuation field with respect to initial data of the height function (see Lemma \ref{lemma:ComparisonMicro}), show the approximations of the initial heights propagate with probability 1 even at the level of microscopic height functions; this is what lets us compare non-stationary heights to stationary heights. As the approximation scale $r$ is arbitrary, we get Theorem \ref{theorem:ASHEZRP} and Theorem \ref{theorem:SBEZRP} under possibly random initial data. This covering argument is rather standard provided our proofs of Theorem \ref{theorem:ASHEZRP} and Theorem \ref{theorem:SBEZRP}, so we do not include it in this paper, but we make this remark in case of possible interest.

We now provide background and history concerning the monotonicity property of the height function process. This ultimately falls out of the attractiveness of the zero-range process, and attractiveness of interacting particle systems has been used crucially in previous literature as a microscopic property of an interacting particle system that reflects an important property of a scaling limit. 
\begin{itemize}[leftmargin=*]
\item Let us take TASEP, which is an interacting particle system that has an attractiveness-type property. In \cite{F1}, the hydrodynamic limit of TASEP, under hyperbolic Eulerian scaling, is shown to be a nonlinear Euler equation from fluid mechanics. Tightness of the space-time particle density is rather straightforward. To identify the limit points as weak solutions to the nonlinear Euler equation, the main obstruction is a lack of uniqueness of solutions to the nonlinear Euler equation unless one restricts to a class of weak solutions satisfying the entropy condition. The attractiveness of TASEP fortunately corresponds to a microscopic version of this entropy condition, which allows \cite{F1} to succeed. We refer to \cite{F2} for a similar implementation of attractiveness of TASEP but for fluctuations of the hydrodynamic limit rather than the hydrodynamic limit itself.
\end{itemize}
\subsection{Non-Simple Exclusion Processes}
We conclude this introduction with an explicit interacting particle system for which our methods apply; by method, we mean comparing general deterministic initial data to almost stationary data that one can then derive density fluctuations for with ideas of energy solution theory \cite{GJS15} or standard extensions thereof. This is to emphasize the generality of the method developed herein. We start by highlighting three key properties of zero-range processes considered herein that allow our method to work. The first is an explicit invariant measure that is sufficiently mixing, so the initial height function is a Brownian bridge on $\mathbb{T}$. The second is monotonicity of height functions, so if two height functions are initially uniformly close, then there is a coupling that keeps the height functions uniformly close (see Lemma \ref{lemma:ComparisonMicro}), and this is where attractivity is useful. The first Brownian invariant measure is also useful in constructing approximations of initial height functions with respect to measures very close to invariant measures, and monotonicity allows us to globally propagate comparison. The third key property is a ``gradient condition", which is necessary even for stationary models using energy solution theory \cite{GJ15, GJ16, GJS15}. We note, however, that this gradient condition is not directly used in this paper, namely for comparison between general and almost stationary data; it is only then used for calculations already done in \cite{GJS15} for stationary models and thus will not appear in our analysis or discussion.

Let us now consider the following non-simple exclusion processes. Fix a maximal jump-length $\mathfrak{m}\geq1$, and take any configuration of particles on $\mathbb{T}_{N}$ with at most one particle per point. The symmetric part of the process is given by, at speed $N^{2}$, choosing a direction (left or right) uniformly at random with probability $1/2$ each, and then jumping to the closest empty point in that direction within distance $\mathfrak{m}$; if no such point exists, then no jump is taken. The asymmetric part of the process is to jump to the closest empty point to the right within distance $\mathfrak{m}$ at speed $N^{3/2}$ (only if such a point exists). The generator of this process is given by, for $\eta\in\{0,1\}^{\mathbb{T}_{N}}$ with $\eta_{x}=1$ meaning a particle at $x$, the formula
\begin{align}
\mathscr{L}^{\mathrm{exc}} \ = \ 2^{-1}N^{2}\sum_{k=1}^{\mathfrak{m}}\sum_{x\in\mathbb{T}_{N}}\prod_{j=1}^{k-1}\eta_{x+j}\mathscr{L}^{\mathrm{exc}}_{x,k} \ + \ N^{\frac32}\sum_{k=1}^{\mathfrak{m}}\sum_{x\in\mathbb{T}_{N}}\eta_{x}(1-\eta_{x+k})\prod_{j=1}^{k-1}\eta_{x+j}\mathscr{L}^{\mathrm{exc}}_{x,k},
\end{align}
where $\mathscr{L}^{\mathrm{exc}}_{x,k}$ is the speed-1 generator for symmetric exclusion on the bond $\{x,x+k\}$. It is easy to check the gradient condition for $\mathscr{L}^{\mathrm{exc}}$, namely $\prod_{j=1}^{k-1}\eta_{x+j}\cdot(\eta_{x+k}-\eta_{x})$ is a gradient. As in \cite{GJ15}, this implies that the product Bernoulli measures of constant density are invariant. Lastly, this process can be checked to satisfy a monotonicity in height functions with a coupling between two configurations given by particles at the same site jumping together whenever possible (of possibly different jump lengths). Indeed, the height function process consists of swapping peaks/valleys, thus monotonicity follows similar to how it does in ASEP in \cite{QS}. 
\subsection{Notation}
First, given any $a,b \in \R$, define $\llbracket a,b \rrbracket \overset{\bullet}= [a,b] \cap \Z$. Second, given any length $k \in \Z$, we define the following discrete gradient operators for $X \in \mathbb{T}_{N}$ and $\varphi: \mathbb{T}_{N} \to \R$:
\begin{align}
\grad_{k} \varphi(X) \ \overset{\bullet}= \ \varphi(X+k) - \varphi(X); \quad \grad_{k}^{N} \ \overset{\bullet}= \ N \grad_{k}.
\end{align}
Deterministic and near-stationary measures we operate with will be denoted by $\{\mathbf{P}_{0,N}\}_{N\in\Z_{>0}}$, and the density fluctuation fields with initial condition sampled via $\{\mathbf{P}_{0,N}\}_{N\in\Z_{>0}}$ will be denoted by $\{\mathbf{Y}^{N,0}\}_{N\in\Z_{>0}}$. In similar spirit, the height functions with initial condition sampled via $\{\mathbf{P}_{0,N}\}_{N\in\Z_{>0}}$ will be denoted by $\{\mathbf{h}^{N,0}\}_{N\in\Z_{>0}}$.
\subsection{Acknowledgements}
The author thanks Amir Dembo and Li-Cheng Tsai and for helpful conversations. The author also thanks Alexander Dunlap and Jimmy He for interesting discussions.
%
%
%
\section{Almost-Stationary Approximations}
The first ingredient towards SBE fluctuations of the density field consists of approximating the associated height function $\mathbf{h}^{N,0}$ by those arising from almost-stationary particle systems in the sense of bounded relative entropy; see Theorem \ref{theorem:ASHEZRPEq}, for example.
\begin{itemize}[leftmargin=*]
\item First assume that the total particle density associated to the initial non-equilibrium measure $\mathbf{P}_{N,0}$ is equal to $\sigma$ for simplicity. We address the necessary changes at the end of this outline, which are basically remedied by an elementary perturbation.

\item The grand-canonical ensemble $\mathbf{P}_{\alpha_{\sigma}}$ on $\Omega_{\mathbb{T}_{N}}$ gives an invariant measure for the particle system with particle density equal to $\sigma$ and with fluctuations on the order $|\mathbb{T}_{N}|^{1/2}$. The order of fluctuations follows by the classical central limit theorem and iid property of one-point marginals of $\mathbf{P}_{\alpha_{\sigma}}$. Preliminarily, we will first condition for these fluctuations to be small with respect to the aforementioned scale as to agree with the particle density with respect to the non-equilibrium initial condition $\mathbf{P}_{N,0}$. In particular, provided any \emph{fixed} $\e \in \R_{>0}$, we construct the new probability measure $\mathbf{P}_{N,\sigma,1,\e}$ on $\Omega_{\mathbb{T}_{N}}$ defined by conditioning the probability measure $\mathbf{P}_{\alpha_{\sigma}}$ on the event in which fluctuations of the particle number are on the order of $\e\sqrt{|\mathbb{T}_{N}|}$, for example. Lemma \ref{lemma:ASA1} below shows that $\mathbf{P}_{N,\sigma,1,\e}$ is stable in the large-$N$ limit with respect to $\mathbf{P}_{\alpha_{\sigma}}$ in the sense of relative entropy. Thus, we may now inherit all of the convergence results in \cite{GJS15} from $\mathbf{P}_{\alpha_{\sigma}}$ to $\mathbf{P}_{N,\sigma,1,\e}$.

\item Under the new probability measure $\mathbf{P}_{N,\sigma,1,\e}$ on $\Omega_{\mathbb{T}_{N}}$ constructed in the previous bullet point, the height function, denoted by $\mathbf{h}^{N,\e}$, under $\mathbf{P}_{N,\sigma,1,\e}$ converges in functional CLT sense to a Brownian motion on $[0,1]$ of $\sigma$-dependent diffusion constant \emph{conditioned} with endpoint in the interval $[-\e,\e] \subseteq \R$. This second bullet point is purely to emphasize this.

\item We proceed to describe the almost-stationary approximations. The height functions $\mathbf{h}^{N,0}$ under the non-equilibrium initial condition $\mathbf{P}_{N,0}$ provide continuous functions on the interval $[0,1] \subseteq \R_{\geq 0}$ upon linear interpolation. By assumption on the particle density under $\mathbf{P}_{N,0}$, we have $\mathbf{h}_{0,0}^{N,0} = \mathbf{h}_{0,|\mathbb{T}_{N}|}^{N,0} = 0$. Courtesy of classical stochastic calculus, an elementary calculation shows that a Brownian motion indexed by $[0,1]$ conditioned to have its endpoint contained in $[-\e,\e] \subseteq \R$ has probability uniformly bounded below by $\kappa_{\e} \in \R_{>0}$ of staying within $\e \in \R_{>0}$ of any pre-determined continuous function $\varphi$ on $[0,1]$ that satisfies the two deterministic constraints $\varphi(0) = 0$ and $|\varphi(1)|\leq\e$ that we condition the Brownian motion on.

\item Given the previous pair of observations, our almost-stationary approximations primarily consist of conditioning $\mathbf{P}_{N,\sigma,1,\e}$ on the event in which the associated height functions $\mathbf{h}^{N,\e}$ are uniformly within $\e \in \R_{>0}$ of $\mathbf{h}^{N,0}$. Fortunately, this new measure is stable in the sense of relative entropy to $\mathbf{P}_{N,\sigma,1,\e}$ as shown in Lemma \ref{lemma:ASA2} below.
\item We now discuss removing the assumption of $\mathbf{h}^{N,0}_{0,|\mathbb{T}_{N}|} = 0$. Consider $\mathbf{h}_{0,|\mathbb{T}_{N}|}^{N,0} = \mathbf{h}_{N} \in \R$, so that $\mathbf{h}_{N}$ is uniformly bounded in $N \in \Z_{>0}$  with a limit $\lim_{N \to \infty} \mathbf{h}_{N} \to \mathbf{h}_{\infty} = \mathbf{h}_{0,1}^{\infty} \in \R$ which holds by assumed convergence $\mathbf{h}^{N,0}\to\mathbf{h}^{\infty}$. For this scenario, we instead construct the probability measure $\mathbf{P}_{N,\sigma,1,\e}$ by conditioning on the fluctuations of the particle number to be on scale $\e|\mathbb{T}_{N}|^{-\frac12}$ \emph{around} $\mathbf{h}_{\infty} \in \R$. The resulting height function $\mathbf{h}^{N,\e}$ thereby converges, again in the sense of functional CLT, to Brownian motion of identical $\sigma$-dependent variance conditioned to have endpoint in $\mathbf{h}_{\infty} + [-\e,\e] \subseteq \R$. The almost-stationary approximations are then constructed in identical fashion by further conditioning on uniform $\e$-closeness to $\mathbf{h}^{N,0}$.
\end{itemize}
We emphasize particle dynamics themselves are entirely irrelevant in the entirety of the section. In particular, all calculations are with respect to the grand-canonical ensembles and hold equally well for the grand-canonical ensembles given by product Bernoulli measures in the case of ASEP, for example.
\subsection{The Rigorous Construction}
We first construct a one-parameter family of events and of measures $\{\mathbf{P}_{N,\sigma,1,\e}\}_{\e\in\R_{>0}}$.
\begin{definition}
Provided any $\e \in \R_{>0}$, recall $\mathbf{h}_{\infty} = \mathbf{h}^{\infty}_{0,1} = \lim_{N\to\infty}\mathbf{h}_{0,|\mathbb{T}_{N}|}^{N,0}$ and define the event
\begin{align}
\mathscr{D}_{N,\e} \ \overset{\bullet}= \ \left\{ \eta \in \Omega_{\mathbb{T}_{N}}: \ \left(|\mathbb{T}_{N}|^{-\frac12} {\sum}_{x \in \mathbb{T}_{N}} \eta_{x}\right) - \sigma |\mathbb{T}_{N}|^{\frac12} \in [\mathbf{h}_{\infty}-\e,\mathbf{h}_{\infty}+\e] \right\}.
\end{align}
We recall $|\mathbb{T}_{N}| =  N$. Moreover, provided any $\e \in \R_{>0}$, define the probability measure
\begin{align}
\d\mathbf{P}_{N,\sigma,1,\e} \ \overset{\bullet}= \ \mathscr{Z}_{\e,N,1}^{-1} \mathbf{1}_{\mathscr{D}_{N,\e}} \d\mathbf{P}_{\alpha_{\sigma}}, \quad \mathscr{Z}_{\e,N,1} \ \overset{\bullet}= \ \mathbf{P}_{\alpha_{\sigma}}\left( \mathscr{D}_{N,\e} \right).
\end{align}
\end{definition}
The first result for the current subsection justifies our consideration of this family $\{\mathbf{P}_{N,\sigma,1,\e}\}_{\e \in \R}$ of probability measures on $\Omega_{\mathbb{T}_{N}}$. Roughly, it allows us to obtain the desired convergence results from \cite{GJS15} though for the family $\{\mathbf{P}_{N,\sigma,1,\e}\}_{\e \in \R}$ rather than the grand-canonical ensemble $\mathbf{P}_{\alpha_{\sigma}}$ by means of a relative entropy estimate and the last statements of Theorems \ref{theorem:ASHEZRPEq} and \ref{theorem:SBEZRPEq}.
\begin{lemma}\label{lemma:ASA1}
Provided any $\e \in \R_{>0}$, we have the following upper bound on relative entropy with implied constant depending only on $\e \in \R_{>0}$ and, in particular, independent of $N \in \Z_{>0}$; see Definition \ref{definition:RE} for the definition of relative entropy:
\begin{align}
\mathsf{H}_{\mathbf{P}_{\alpha_{\sigma}}}(\mathbf{P}_{N,\sigma,1,\e}) \ &\lesssim_{\e} \ 1.
\end{align}
\end{lemma}
\begin{proof}
By definition of relative entropy, we first have the following preliminary bound with universal implied constant:
\begin{align}
\mathsf{H}_{\mathbf{P}_{\alpha_{\sigma}}}(\mathbf{P}_{N,\sigma,1,\e})  \ &\lesssim \ 1 \ + \ \E^{\mathbf{P}_{N,\sigma,1,\e}} \mathbf{1}_{\mathscr{D}_{N,\e}} \left( \log \frac{\d\mathbf{P}_{N,\sigma,1,\e}}{\d\mathbf{P}_{\alpha_{\sigma}}} \right)_{+}. \label{eq:ASA11}
\end{align}
In \eqref{eq:ASA11}, we have adopted the notation $x_{+} = x \vee 0$ for the non-negative part of a number. Indeed, whenever the logarithmic-factor $\log \varphi$, with $\varphi$ the Radon-Nikodym derivative at hand, in the relative entropy is negative, the entire factor of $\varphi \log \varphi$ is uniformly bounded. Note measure $\mathbf{P}_{N,\sigma,1,\e}$ is supported on the event $\mathscr{D}_{N,\e}$ by definition. Second, by definition of $\mathbf{P}_{N,\sigma,1,\e}$ we get
\begin{align}
\mathbf{1}_{\mathscr{D}_{N,\e}} \left( \log \frac{\d\mathbf{P}_{N,\sigma,1,\e}}{\d\mathbf{P}_{\alpha_{\sigma}}} \right)_{+} \ = \ \mathbf{1}_{\mathscr{D}_{N,\e}}\log \mathscr{Z}_{\e,N,1}^{-1} \ &\leq \ \log \mathscr{Z}_{\e,N,1}^{-1}. \label{eq:ASA12}
\end{align}
Combining \eqref{eq:ASA11} and \eqref{eq:ASA12}, it suffices to prove that $\mathscr{Z}_{\e,N,1} \gtrsim_{\e} 1$ with an implied constant depending only on $\e \in \R_{>0}$ and thus independent of $N \in \Z_{>0}$. To this end, we observe that under the probability measure $\mathbf{P}_{\alpha_{\sigma}}$, the following functional and the $\eta$-variables all admit more than 2 moments per the moment bound in Definition \ref{definition:intro1}, and therefore
\begin{align}
\left(|\mathbb{T}_{N}|^{-\frac12} {\sum}_{x \in \mathbb{T}_{N}} \eta_{x}\right) - \sigma |\mathbb{T}_{N}|^{\frac12}
\end{align}
converges, in law, in the large-$N$ limit to the Gaussian distribution of mean 0 and variance $\mathrm{v}_{\sigma} \in \R_{>0}$ depending only on $\sigma\in\R$. If $\mathscr{N}_{\sigma}$ denotes a random variable with such Gaussian distribution, as $\mathrm{v}_{\sigma}$ is uniformly positive and bounded in $N\in\Z_{>0}$ we have
\begin{align}
\mathscr{Z}_{\e,N,1} \ = \ \mathbf{P}\left(\mathscr{N}_{\sigma} \in [\mathbf{h}_{\infty}-\e,\mathbf{h}_{\infty}+\e]\right) \ + \ \kappa_{N,\sigma,\mathbf{h}_{\infty},\e}, \label{eq:ASA1CLT}
\end{align}
where $\lim_{N \to \infty} \kappa_{N,\sigma,\mathbf{h}_{\infty},\e} = 0$. The desired lower bound $\mathscr{Z}_{\e,N,1} \gtrsim_{\e} 1$ now follows from an elementary lower bound on the Gaussian probability in \eqref{eq:ASA1CLT} and afterwards choosing $N \in \Z_{>0}$ sufficiently large to make $\kappa_{N,\sigma,\mathbf{h}_{\infty},\e}$ sufficiently small.
\end{proof}
We now introduce a further constraint on the allowed configurations in $\Omega_{\mathbb{T}_{N}}$ on which $\mathbf{P}_{N,\sigma,1,\e}$ is supported.
\begin{definition}
Provided any $\e \in \R_{>0}$, define the event
\begin{align}
\mathscr{E}_{N,\e} \ &\overset{\bullet}= \ \left\{ \eta \in \Omega_{\mathbb{T}_{N}}: \ \|\mathbf{h}_{\eta,x} - \mathbf{h}_{0,N^{-1}x}^{\infty}\|_{\mathscr{L}^{\infty}_{x\in\mathbb{T}_{N}}} \ \leq \ \e \right\},
\end{align}
in which we have defined the following function $\mathbf{h}_{\eta,\bullet}:\mathbb{T}_{N} \to \R$ which extends to $[1,N] \to \R$ by linear interpolation:
\begin{align}
\mathbf{h}_{\eta,x} \ \overset{\bullet}= \ N^{-\frac12} {\sum}_{1\leq y \leq x} \left( \eta_{y} - \sigma \right).
\end{align}
Moreover, define the probability measure
\begin{align}
\d\mathbf{P}_{N,\sigma,2,\e} \ \overset{\bullet}= \ \mathscr{Z}_{\e,N,2}^{-1} \mathbf{1}_{\mathscr{E}_{N,\e}} \d\mathbf{P}_{N,\sigma,1,\e}, \quad \mathscr{Z}_{\e,N,2} \ \overset{\bullet}= \ \mathbf{P}_{N,\sigma,1,\e}\left(\mathscr{E}_{N,\e}\right).
\end{align}
\end{definition}
\noindent We show $\mathbf{P}_{N,\sigma,2,\e}$ is stable with respect to $\mathbf{P}_{N,\sigma,1,\e}$ in relative entropy, so convergences in \cite{GJS15} hold for initial measure $\mathbf{P}_{N,\sigma,2,\e}$.
\begin{lemma}\label{lemma:ASA2}
Provided any $\e \in \R_{>0}$, we have the relative entropy bound $\mathsf{H}_{\mathbf{P}_{N,\sigma,1,\e}}(\mathbf{P}_{N,\sigma,2,\e}) \lesssim_{\e} 1$.
\end{lemma}
\begin{proof}
Exactly as in the proof of Lemma \ref{lemma:ASA1}, it suffices to show $\mathscr{Z}_{\e,N,2} \gtrsim_{\e} 1$. We first make the following observation.
\begin{itemize}[leftmargin=*]
\item Under $\mathbf{P}_{N,\sigma,1,\e}$, the process $\{\mathbf{h}_{\eta,x}\}_{x\in\mathbb{T}_{N}}$ is a discrete martingale with index $x\in\mathbb{T}_{N}$ with respect to the canonical filtration with iid individual step distributions having $2+\e_{0}$-moments then \emph{conditioned} to have its endpoint in $[\mathbf{h}_{\infty}-\e,\mathbf{h}_{\infty}+\e]$. The step distributions are each scaled by $N^{-\frac12} \asymp |\mathbb{T}_{N}|^{-\frac12}$.
\end{itemize}
Combining this observation with a moment calculation exactly as in the proof of Lemma \ref{lemma:ASA1}, we obtain that under $\mathbf{P}_{N,\sigma,1,\e}$, the process $\{\mathbf{h}_{\eta,Nx}\}_{x\in\llbracket0,1\rrbracket}$ converges uniformly to Brownian motion of $\sigma$-dependent diffusion coefficient in the large-$N$ limit \emph{conditioning} on its endpoint to be in $[\mathbf{h}_{\infty}-\e,\mathbf{h}_{\infty}+\e]$. We denote by $\mathfrak{B}$ a Brownian-type stochastic process of this law. Next, for $\varrho \in \R_{>0}$, it will be convenient to consider $\mathbb{T}_{N,\e,\varrho} \overset{\bullet}= N\e^{\varrho}\Z \cap [0,N]$.

As a final preliminary, let us define the following auxiliary events provided any sufficiently large but universal constants $\varrho,\mathscr{C} \in \R_{>0}$ and sufficiently small but universal constant $\nu\in\R_{>0}$, all in a manner to be made precise shortly. 
\begin{subequations}
\begin{align}
\mathscr{E}_{N,\e,1,\varrho} \ &\overset{\bullet}= \ \left\{ {\sup}_{x\in\mathbb{T}_{N,\e,\varrho}}|\mathfrak{B}_{N^{-1}x} - \mathbf{h}_{0,N^{-1}x}^{\infty}| \ \leq \ \frac13\e\right\}; \\
\mathscr{E}_{N,\e,2,\varrho} \ &\overset{\bullet}= \ \left\{ \eta \in \Omega_{\mathbb{T}_{N}}: \ {\sup}_{x\in\mathbb{T}_{N,\e,\varrho}}|\mathbf{h}_{\eta,x} - \mathbf{h}_{0,N^{-1}x}^{\infty}| \ \leq \ \frac13 \e \right\} ;\\
\mathscr{E}_{N,\e,3,\mathscr{C}} \ &\overset{\bullet}= \ \left\{ \eta \in \Omega_{\mathbb{T}_{N}}: \ {\sup}_{x\neq y\in\mathbb{T}_{N}}|\mathbf{h}_{\eta,x}-\mathbf{h}_{\eta,y}| \cdot |x-y|^{-\nu} \ \leq \ \mathscr{C} N^{-\frac12\nu} \right\}.
\end{align}
\end{subequations}
Observe the event $\mathscr{E}_{N,\e,1,\varrho}$ is unambiguously defined once we have specified a probability measure for the Brownian motion $\mathfrak{B}$ because $\mathbf{h}^{\infty}$ is deterministic. The same applies to $\mathscr{E}_{N,\e,2,\varrho}$. The parameter $\nu \in \R_{>0}$ in $\mathscr{E}_{N,\e,3,\varrho}$ is a universal Holder exponent.

First, we observe the following lower bound for $\mathscr{Z}_{\e,N,2}$ provided any $\mathscr{C} \in \R_{>0}$ sufficiently large but universal and $N \in \Z_{>0}$ sufficiently large depending only on $\nu,\e\in\R_{>0}$ and $\varrho \in \R_{>0}$ sufficiently large but universal depending only on $\nu \in \R_{>0}$:
\begin{align}
\mathscr{Z}_{\e,N,2} \ &\geq \ \mathbf{P}_{N,\sigma,1,\e} \left(\mathscr{E}_{N,\e,2,\varrho} \cap \mathscr{E}_{N,\e,3,\mathscr{C}} \right) \\
&= \ \mathbf{P}_{N,\sigma,1,\e} \left(\mathscr{E}_{N,\e,2,\varrho} \right] \ - \ \mathbf{P}_{N,\sigma,1,\e} \left(\mathscr{E}_{N,\e,2,\varrho} \cap \mathscr{E}_{N,\e,3,\mathscr{C}}^{C} \right) \\
&\geq \ \mathbf{P}_{N,\sigma,1,\e} \left(\mathscr{E}_{N,\e,2,\varrho} \right) \ - \ \mathbf{P}_{N,\sigma,1,\e} \left(\mathscr{E}_{N,\e,3,\mathscr{C}}^{C} \right). \label{eq:ASA21}
\end{align}
The previous lower bound implicitly uses the \emph{deterministic} regularity of $\mathbf{h}^{\infty}$, which is random at the microscopic level for $\mathbf{h}_{\eta,\bullet}$ as indicated by the necessity of the event $\mathscr{E}_{N,\e,3,\mathscr{C}}$. In words, this lower bound is saying that the difference between $\mathbf{h}_{\eta,\bullet}$ and $\mathbf{h}^{\infty}_{0,N^{-1}\bullet}$ on $\mathbb{T}_{N}$ is controlled by the difference between the two on the scale $\e^{\varrho}$-lattice $\mathbb{T}_{N,\e,\varrho}$ and because regularity of $\mathbf{h}^{\infty}$ and $\mathbf{h}_{\eta,\bullet}$ let us bootstrap from the scale $N\e^{\varrho}$-lattice $\mathbb{T}_{N,\e,\varrho}$ to all of $\mathbb{T}_{N}$.

Recall the uniform convergence $\mathbf{h}_{\eta,Nx} \to \mathfrak{B}_{x}$ in law in the large-$N$ limit. Thus, the $\mathscr{E}_{N,\e,2,\varrho}$-probability is asymptotically in the large-$N$ limit equal to the probability of $\mathscr{E}_{N,\e,1,\varrho}$ as $\mathscr{E}_{N,\e,1,\varrho}$ and $\mathscr{E}_{N,\e,2,\varrho}$ look at a scale $\e^{\varrho}$ discretization of $[0,1]$ that is independent of $N\in\Z_{>0}$; we note that we would not necessarily have convergence of probabilities for events that depend on the scaling parameter $N$ under which $\mathbf{h}_{\eta,Nx}$ converges to the Brownian motion $\mathfrak{B}_{x}$:
\begin{align}
\mathbf{P}_{N,\sigma,1,\e} \left( \mathscr{E}_{N,\e,2,\varrho} \right) \ &= \ \mathbf{P}_{N,\sigma,1,\e} \left( \mathscr{E}_{N,\e,1,\varrho} \right) \ + \ \kappa_{N,\sigma,\mathbf{h}_{\infty},\e,\varrho} \\
&= \ C(\e,\varrho) \ - \ |\kappa_{N,\sigma,\mathbf{h}_{\infty},\e,\varrho}|, \label{eq:ASA22}
\end{align}
where $\lim_{N\to\infty}\kappa_{N,\sigma,\mathbf{h}_{\infty},\e,\varrho} = 0$ given any $\varrho \in \R_{>0}$ universal. The lower bound \eqref{eq:ASA22} follows by a stochastic calculus support lemma in Lemma \ref{lemma:SupportTheorem}. For $N\gtrsim_{\e}1$, we may drop the $|\kappa_{N,\sigma,\mathbf{h}_{\infty},\e,\varrho}|$-term in \eqref{eq:ASA22} and update $C(\e,\varrho)\to\frac12C(\e,\varrho)$.

Combining \eqref{eq:ASA21} and \eqref{eq:ASA22}, it suffices to get $\mathbf{P}_{N,\sigma,1,\e} [\mathscr{E}_{N,\e,3,\mathscr{C}}^{C}] \to_{\mathscr{C}\to\infty} 0$ uniformly in $N \in \Z_{>0}$ and for any given $\e \in \R_{>0}$ fixed. We emphasize that this amounts to controlling regularity of a linearly interpolated $\mathbf{h}_{\eta,\bullet}$ at a microscopic level under the $\mathbf{P}_{N,\sigma,1,\e}$-measure. By the entropy inequality in Appendix 1.8 of \cite{KL} and Lemma \ref{lemma:ASA1}, it suffices to get this probability vanishes in the large-$\mathscr{C}$ limit uniformly in $N\in\Z_{>0}$ for any $\e\in\R_{>0}$ fixed but now under the grand-canonical ensemble $\mathbf{P}_{\alpha_{\sigma}}$ instead of under $\mathbf{P}_{N,\sigma,1,\e}$. Under the grand-canonical ensemble, the height function $\mathbf{h}_{\eta,\bullet}$ is actually a martingale free of any conditioning with index $\bullet\in\mathbb{T}_{N}$ with iid increments of uniformly bounded $2+\delta$-moment for $\delta\in\R_{>0}$ universal courtesy of Definition \ref{definition:intro1}, then scaled by $N^{-\frac12} \asymp |\mathbb{T}_{N}|^{-\frac12}$, so controlling regularity is standard via a Kolmogorov continuity criterion.
\end{proof}
As an immediate consequence of Lemma \ref{lemma:ASA1} and Lemma \ref{lemma:ASA2}, we establish the following result which gives the final input from the current section that we will use towards the proof of Theorem \ref{theorem:SBEZRP}.
\begin{corollary}\label{corollary:ASA}
Provided any $\e \in \R_{>0}$, there exists a probability measure $\mathbf{P}_{N,\e}$ on $\Omega_{\mathbb{T}_{N}}$ which satisfies the following constraints.
\begin{itemize}[leftmargin=*]
\item \emph{Theorem \ref{theorem:SBEZRPEq}} holds upon replacing $\mathbf{P}_{\alpha_{\sigma}}$ by $\mathbf{P}_{N,\e}$.
\item Outside an event of probability $\delta_{N} \to_{N \to \infty} 0$ under $\mathbf{P}_{0,N}^{\e,\otimes} \overset{\bullet}= \mathbf{P}_{0,N} \otimes\mathbf{P}_{N,\e}$, with universal implied constant we have
\begin{align}
\| \mathbf{h}_{0,x}^{N,0} - \mathbf{h}_{0,x}^{N,\e} \|_{\mathscr{L}^{\infty}_{x\in\mathbb{T}_{N}}} \ &\lesssim \ \e,
\end{align}
where $\mathbf{h}^{N,0}$ is the height function sampled from $\eta_{1} \sim \mathbf{P}_{0,N}$ and $\mathbf{h}^{N,\e}$ is the height function sampled from $\eta_{2} \sim \mathbf{P}_{N,\e}$.
\end{itemize}
\end{corollary}
\begin{proof}
Consider the measure $\mathbf{P}_{N,\e} = \mathbf{P}_{N,\sigma,2,\e}$ provided any $\e \in \R_{>0}$.
\begin{itemize}[leftmargin=*]
\item By the relative entropy bound in Lemma \ref{lemma:ASA1}, Theorem \ref{theorem:SBEZRPEq} holds upon replacing $\mathbf{P}_{\alpha_{\sigma}}$ with $\mathbf{P}_{N,\sigma,1,\e}$. By the relative entropy estimate in Lemma \ref{lemma:ASA2} and identical reasoning, we deduce that Theorem \ref{theorem:SBEZRPEq} holds upon replacing $\mathbf{P}_{\alpha_{\sigma}}$ with $\mathbf{P}_{N,\sigma,2,\e}$.
\item It remains to establish the necessary estimate on height functions. First, we have the following for $\delta_{N} \to_{N \to \infty} 0$ simply by what it means to have deterministic and near-stationary initial data for the particle system and height function:
\begin{align}
\mathbf{P}_{0,N}\left(\| \mathbf{h}_{0,x}^{N,0} - \mathbf{h}_{0,N^{-1}x}^{\infty} \|_{\mathscr{L}^{\infty}_{x\in\mathbb{T}_{N}}} \ \leq \ \e\right) \ \geq \ 1 \ - \ \delta_{N}
\end{align}
Moreover, by definition of the auxiliary measure $\mathbf{P}_{N,\sigma,2,\e}$ through the event $\mathscr{E}_{N,\e}$, an identical probability-1 statement holds for $\mathbf{P}_{N,\sigma,2,\e}$-height functions. Precisely, if $\mathbf{h}^{N,\e}$ denotes the height function sampled according to $\mathbf{P}_{N,\sigma,2,\e}$, we have
\begin{align}
\mathbf{P}_{N,\sigma,2,\e}\left(\| \mathbf{h}_{0,x}^{N,\e} - \mathbf{h}_{0,N^{-1}x}^{\infty} \|_{\mathscr{L}^{\infty}_{x\in\mathbb{T}_{N}}} \ \leq \ \e\right) \ = \ 1.
\end{align}
As $\mathbf{h}^{\infty}$ is deterministic, the previous statements and the triangle inequality give the desired height function bound.
\end{itemize}
This completes the proof.
\end{proof}
%
%
%
\section{Comparison Principles}
For the present section, we exploit both the comparison principle for the SHE from \cite{Mu} and a suitable microscopic version of this principle at the microscopic level of the particle system, though in reverse order of what we just stated.
\subsection{Microscopic Principle}
The microscopic analog of the comparison principle for the SHE that we introduce here serves to propagate the height function bound in Corollary \ref{corollary:ASA} globally in time with respect to a coupling of identical ZRP-dynamics. Keeping in mind that the SHE comparison principle in \cite{Mu} couples two solutions to SHE to the same noise, we employ the \emph{basic coupling} that couples two different initial configurations through Poisson clocks; see Section 1.4 in \cite{BGS}. Roughly speaking, the basic coupling above between any particle configurations attaches the same Poisson clock to each directed bond in the lattice. Particles jump together between the two configurations whenever they are able to. Attractiveness of the particle system implies that if one particle configuration observes a particle-jump while the other does not, then before this jump, the particle number at the site from which the particle jumped is strictly greater in the former configuration than the latter. We note that after this jump, the number of particles at the site from where the particle jumped is at least that in the configuration that did not see a jump. Indeed, the number was strictly greater before the jump, and the jump only removed one particle from the site in the configuration that saw the jump. We refer to \cite{BGS} for more background on this basic coupling.
\begin{remark}\label{remark:InfCoupling}
We couple all $\e \in \R_{\geq 0}$-measures simultaneously on the same probability space.
\end{remark}
The following result asserts that under the above basic coupling, the height function estimate in Corollary \ref{corollary:ASA} propagates globally in time, providing a genuine analog of the comparison principle in \cite{Mu}.
\begin{lemma}\label{lemma:ComparisonMicro}
\fontdimen2\font=1.7pt Consider any arbitrarily small although universal parameter $\e \in \R_{>0}$. Moreover, we define $\mathbf{P}_{T,N}^{\e,\otimes}$ as the probability measure on $\mathscr{D}([0,\R_{\geq0}],\Omega_{\mathbb{T}_{N}}^{2})$ defined by the stochastic process associated to the basic coupling dynamic with initial condition the product measure $\mathbf{P}_{0,N}^{\e,\otimes}$ from \emph{Corollary \ref{corollary:ASA}}. 

With probability at least $1 - \delta_{N}$, where $\delta_{N} \in \R_{>0}$ is the same constant from \emph{Corollary \ref{corollary:ASA}}, with respect to $\mathbf{P}_{T,N}^{\e,\otimes}$, we have the following estimate for height functions with the same universal implied constant as in \emph{Corollary \ref{corollary:ASA}}:
\begin{align}
\| \mathbf{h}_{S,x}^{N,0} - \mathbf{h}_{S,x}^{N,\e} \|_{\mathscr{L}^{\infty}_{S\geq0,x\in\mathbb{T}_{N}}} \ &\lesssim \ \e.
\end{align}
Above $\mathbf{h}^{N,0}$ is the height function associated to the projection of the basic coupling onto the first-coordinate with initial condition sampled from $\mathbf{P}_{0,N}$, and $\mathbf{h}^{N,\e}$ is the height function associated to the projection onto the second coordinate with initial condition sampled from the probability measure $\mathbf{P}_{N,\e}$ within \emph{Corollary \ref{corollary:ASA}}.
\end{lemma}
\begin{proof}
By Corollary \ref{corollary:ASA}, we first get the following two-sided estimate with universal constant $\kappa' \in \R_{>0}$ uniformly in $x \in \mathbb{T}_{N}$ and for $T = 0$ with probability at least $1 - \delta_{N}$ under $\mathbf{P}_{0,N}^{\e,\otimes}$:
\begin{align}
\mathbf{h}_{T,x}^{N,\e} - \kappa' \e \ \leq \ \mathbf{h}_{T,x}^{N,0} \ \leq \ \mathbf{h}_{T,x}^{N,\e} + \kappa' \e.
\end{align}
We observe that Lemma \ref{lemma:ComparisonMicro} would follow from showing the previous two-sided bound extends globally in time with the same positive constant $\kappa' $ with probability 1 with respect to the randomness of the basic coupling dynamic. Indeed, $\kappa'$ is \emph{independent} of the approximation parameter $\e \in \R_{>0}$ and \emph{independent} of $N\in\Z_{>0}$. Propagation of the above estimate is proven as follows:
\begin{itemize}[leftmargin=*]
\item Consider the first time $T \geq 0$ at which the inequality fails. At this time, there must exist a particle-jump occurring between the sites $x,x+1 \in \mathbb{T}_{N}$ in the $0$-configuration but not the $\e$-configuration, or vice versa.
\item Suppose a particle jumps from $x \rightsquigarrow x+1$ in the $0$-configuration but not the $\e$-configuration. Note $\mathbf{h}_{T,x}^{N,0}$ exhibits a negative-jump due to this particle-action but $\mathbf{h}_{T,x}^{N,\e}$ does not change. Thus, we deduce the upper bound remains valid after this jump. To show the lower bound also remains valid, as noted in the first paragraph of this subsection, the number of particles at site $x \in \mathbb{T}_{N}$ within the $0$-configuration is at least that in the $\e$-configuration. Moreover, because Poisson clocks ring at distinct times with probability 1, the inequality remains valid both before and immediately after this jump at the site $x - 1 \in \mathbb{T}_{N}$; before and after this jump, nothing has happened to left of $x\in\mathbb{T}_{N}$. We deduce the lower bound must remain valid immediately after the jump at site $x \in \mathbb{T}_{N}$ just by definition of the height function, simply adding the $\mathbf{h}^{N,0}$ and $\mathbf{h}^{N,\e}$ values at $x-1\in\mathbb{T}_{N}$ with the number of particles at $x\in\mathbb{T}_{N}$ in the respective configurations, which we compared at the beginning of this paragraph.
\item Similarly, suppose a particle jumps from $x \rightsquigarrow x+1$ in the $\e$-configuration but not the $0$-configuration. In this case, the height function $\mathbf{h}_{T,x}^{N,\e}$ exhibits a negative-jump due to such a particle action, but $\mathbf{h}_{T,x}^{N,0}$ does not change. This implies the lower bound remains valid immediately after this jump. To show the upper bound remains valid, we proceed as in the previous bullet point. The inequality must remain valid at site $x-1 \in \mathbb{T}_{N}$ both before and after the aforementioned particle action, because particles strictly to the left of $x\in\mathbb{T}_{N}$ are unaffected by the previous jump. Furthermore, as consequence of the attractiveness of the zero-range particle system under the basic coupling, the total number of particles at the site $x \in \mathbb{T}_{N}$ within the $\e$-configuration must be at least that in the $0$-configuration. We get from this that the upper bound remains valid immediately after the aforementioned jump as well using reasoning identical to that at the end of the previous bullet point.

\item Considering the remaining cases of a particle jumping $x+1\rightsquigarrow x$ within one configuration but not the other, propagation of the inequality immediately after this jump follows by basically the same arguments in the previous two bullet points. In particular, for these jumps, we note the values of $\mathbf{h}^{N,0}_{T,x+1}$ and $\mathbf{h}^{N,\e}_{T,x+1}$ do not change, whereas the number of particles at $x+1$ in the configuration that sees the jump is at least the number of particles at $x+1$ in the configuration that does not see the jump \emph{after} the jump has been executed. To recover values $\mathbf{h}^{N,0}_{T,x}$ and $\mathbf{h}^{N,\e}_{T,x}$, we then subtract from $\mathbf{h}^{N,0}_{T,x+1}$ and $\mathbf{h}^{N,\e}_{T,x+1}$ these particle numbers at $x+1$ up to a factor of $N^{-1/2}$ and a shared additive factor $N^{-1/2}\sigma$.
\end{itemize}
This completes the proof.
\end{proof}
\subsection{Compactness of SHE}
By Theorem 2.13 in \cite{GP}, the limiting solutions to the stochastic Burgers equation obtained by Theorem \ref{theorem:SBEZRPEq} for particle systems with initial measure $\mathbf{P}_{N,\e}$ correspond to the Cole-Hopf solution of the KPZ equation with the initial condition $\mathbf{h}^{\infty,\e}$ which are uniformly within $\e \in \R_{>0}$ of $\mathbf{h}^{\infty}$. Moreover, Lemma \ref{lemma:ComparisonMicro} above provides a mechanism, which we detail later in the proof of Theorem \ref{theorem:SBEZRP}, to approximate the fluctuation field $\mathbf{Y}^{N,0}$ with initial condition sampled via $\mathbf{P}_{0,N}$ by fluctuation fields $\{\mathbf{Y}^{N,\e}\}_{\e\in\R_{>0}}$ with initial condition sampled with the almost-stationary approximations $\mathbf{P}_{N,\e}$, respectively. In particular, a key step in our proof of Theorem \ref{theorem:SBEZRP} will be the following compactness for the SHE. We set $\mathbf{h}^{\infty,0}=\mathbf{h}^{\infty}$ below.
\begin{lemma}\label{lemma:CompactSHE}
\fontdimen2\font=1.7pt Suppose $\mathbf{h}^{\infty}$ denotes the Cole-Hopf solution to the \emph{KPZ equation \eqref{eq:SPDE}} with the initial condition $\mathbf{h}^{\infty}_{0,\bullet}$. Moreover, for any $\e \in \R_{>0}$, let $\mathbf{h}^{\infty,\e}$ denote the Cole-Hopf solution to the \emph{KPZ equation \eqref{eq:SPDE}} with initial condition $\mathbf{h}_{0,\bullet}^{\infty,\e} = \lim_{N \to \infty} \mathbf{h}_{0,N\bullet}^{N,\e}$, where $\mathbf{h}_{0,N\bullet}^{N,\e}$ is the height function sampled via $\mathbf{P}_{N,\e}$ from \emph{Corollary \ref{corollary:ASA}}. For any $\Phi \in \mathscr{C}^{\infty}(\mathbb{T})$ and any $T \in \R_{\geq 0}$, we have the convergence in probability
\begin{align}
\int_{\mathbb{T}} \Phi(X) \cdot \mathbf{h}_{T,X}^{\infty,\e} \ \d X \ \longrightarrow_{\e \to 0^{+}} \ \int_{\mathbb{T}} \Phi(X) \cdot \mathbf{h}_{T,X}^{\infty} \ \d X.
\end{align}
Here, the solutions $\mathbf{h}^{\infty}$ and $\{\mathbf{h}^{\infty,\e}\}_{\e\in\R_{>0}}$ are jointly coupled on the same probability space by sampling independently initial data and coupling all solutions to the same space-time white noise.
\end{lemma}
\begin{remark}
Our proof of Lemma \ref{lemma:CompactSHE} below could be noticeably simplified upon adopting the regularity structures approach of Hairer in \cite{Hai14}, but this argument only works in compact domains. In view of possible future work for interacting particle systems on non-compact sets, the approach via the Cole-Hopf transform appears most robust. In particular, the following proof of Lemma \ref{lemma:CompactSHE} provides a proof where the torus $\mathbb{T}$ is replaced by the full-line $\R$ upon minor modifications if the space of test functions is compactly supported smooth functions. At the end of the day, however, Lemma \ref{lemma:CompactSHE} is an elementary SPDE exercise.
\end{remark}
\begin{proof}
All probabilities and expectations are to be taken with respect to the joint law of the respective initial data and the common space-time white noise. We also let $\langle\varphi,\psi\rangle$ denote the natural pairing of functions on $\mathbb{T}$ given by integrating the product $\varphi\psi$ against Lebesgue measure. Consider any $\delta \in \R_{>0}$ arbitrarily small though universal and define the sequence of events
\begin{align}
\mathscr{E}_{\e,\delta,\Phi} \ \overset{\bullet}= \ \left\{ | \langle\Phi, \mathbf{h}_{T,X}^{\infty,\e}\rangle \ - \ \langle\Phi, \mathbf{h}_{T,X}^{\infty}\rangle| \ \geq \ \delta \right\}.
\end{align}
It suffices to show that $\mathbf{P}[\mathscr{E}_{\e,\delta,\Phi}] \to_{\e\to0^{+}} 0$ provided any $\delta \in \R_{>0}$ independent of $\e \in \R_{>0}$. To this end, we first observe the following containment of events for $\e \in (0,1]$ and arbitrarily large though universal constants $\mathscr{C},\mathscr{C}' \in \R_{>0}$ subject to the constraint that $\mathscr{C}' \gtrsim_{\mathscr{C}} 1$ is sufficiently large depending only on $\mathscr{C}\in\R_{>0}$:
\begin{align}
\mathscr{E}_{\e,\delta,\Phi}^{C} \ &\supseteq \ \mathscr{E}_{\e,\delta,\Phi,\mathscr{C}}^{1} \cap \mathscr{E}_{\e,\delta,\Phi,\mathscr{C}'}^{2}, \label{eq:CompactSHE1}
\end{align}
where the second event below estimates pairing with $\mathbf{z}^{\infty}$ and with $\mathbf{z}^{\infty,\e}$, whereas the first event below provides a priori estimates for $\mathbf{h}^{\infty,\e}$ that translate comparison between $\mathbf{h}^{\infty,\e}$ and $\mathbf{h}^{\infty}$ into comparison between $\mathbf{z}^{\infty,\e}$ and $\mathbf{z}^{\infty}$ as we explain shortly:
\begin{subequations}
\begin{align}
\mathscr{E}_{\e,\delta,\Phi,\mathscr{C}}^{1} \ &\overset{\bullet}= \ \left\{ -\mathscr{C} \ \leq \ \inf_{\e\in[0,1]} \inf_{S\in[0,T]}\inf_{X\in\mathbb{T}} \mathbf{h}_{S,X}^{\infty,\e} \ \leq \ \sup_{\e\in[0,1]}\sup_{S\in[0,T]}\sup_{X\in\mathbb{T}} \mathbf{h}_{S,X}^{\infty,\e} \ \leq \ \mathscr{C} \right\}; \\
\mathscr{E}_{\e,\delta,\Phi,\mathscr{C}'}^{2} \ &\overset{\bullet}= \ \left\{ \int_{\mathbb{T}}|\mathbf{z}_{T,X}^{\infty} - \mathbf{z}_{T,X}^{\infty,\e}| \ \d X \ \leq \ \frac{1}{\mathscr{C}'\|\Phi\|_{\mathscr{L}^{\infty}_{X}}} \delta \right\}.
\end{align}
\end{subequations}
For a deterministic and uniformly bounded constant $C$ determined by the parameters $\alpha,\beta,\zeta$ in Theorem \ref{theorem:SBEZRP}, we have $\mathbf{z}^{\infty} = C\log \mathbf{h}^{\infty}$ and $\mathbf{z}^{\infty,\e} = C\log \mathbf{h}^{\infty,\e}$. Moreover, on $\mathscr{E}_{\e,\delta,\Phi,\mathscr{C}}^{1}$, the value $\e = 0$ corresponds to the solution $\mathbf{h}^{\infty}$ with deterministic initial data. Indeed, \eqref{eq:CompactSHE1} follows by the upcoming calculation for which we condition on the intersection $\mathscr{E}_{\e,\delta,\Phi,\mathscr{C}}^{1} \cap \mathscr{E}_{\e,\delta,\Phi,\mathscr{C}'}^{2}$:
\begin{align}
| \langle\Phi, \mathbf{h}_{T,X}^{\infty,\e}\rangle \ - \ \langle\Phi, \mathbf{h}_{T,X}^{\infty}\rangle|  \ &\lesssim \ \|\Phi\|_{\mathscr{L}^{\infty}_{X}} \int_{\mathbb{T}} |\log \mathbf{z}_{T,X}^{\infty} - \log \mathbf{z}_{T,X}^{\infty,\e}| \ \d X \ \leq \ \|\Phi\|_{\mathscr{L}^{\infty}_{X}} \mathscr{C}' \int_{\mathbb{T}}|\mathbf{z}_{T,X}^{\infty} - \mathbf{z}_{T,X}^{\infty,\e}| \ \d X;
\end{align}
for the preceding calculation, we employed the elementary inequality $|\log X - \log Y| \leq \mathscr{C}' |X-Y|$ for $\mathrm{e}^{-\kappa'\mathscr{C}} \leq X,Y \leq \mathrm{e}^{\kappa'\mathscr{C}}$ if $\mathscr{C}' \gtrsim_{\mathscr{C}} 1$ is sufficiently large depending only on $\mathscr{C} \in \R_{>0}$ and $\kappa'\in\R_{>0}$ is universal. Thus, from \eqref{eq:CompactSHE1} we obtain
\begin{align}
\mathbf{P}\left(\mathscr{E}_{\e,\delta,\Phi}\right) \ &\leq \ \mathbf{P}\left(\left(\mathscr{E}_{\e,\delta,\Phi,\mathscr{C}}^{1}\right)^{C}\right) \ + \ \mathbf{P}\left(\left(\mathscr{E}_{\e,\delta,\Phi,\mathscr{C}'}^{2}\right)^{C}\right). \label{eq:CompactSHE2}
\end{align}
To estimate the first probability on the RHS of \eqref{eq:CompactSHE2}, observe that for $\kappa'\in\R_{>0}$ a universal constant, we have
\begin{align}
\mathscr{E}_{\e,\delta,\Phi,\mathscr{C}}^{1} \ = \ \left\{ \mathrm{e}^{-\kappa'\mathscr{C}} \ \leq \ \inf_{\e\in[0,1]} \inf_{S\in[0,T]}\inf_{X\in\mathbb{T}} \mathbf{z}_{S,X}^{\infty,\e} \ \leq \ \sup_{\e\in[0,1]}\sup_{S\in[0,T]}\sup_{X\in\mathbb{T}} \mathbf{z}_{S,X}^{\infty,\e} \ \leq \ \mathrm{e}^{\kappa'\mathscr{C}} \right\}.
\end{align}
We now use the continuum comparison principle for the SHE in \cite{Mu}. In particular, observe that with probability 1 with respect to sampling jointly independently for $\e \in [0,1]$, we get the following two-sided bound for $T = 0$ and all $X \in \mathbb{T}$ and $\e \in [0,1]$ simultaneously which follows by Corollary \ref{corollary:ASA} and then taking the large-$N$ limit to deduce that $\mathbf{h}^{\infty,\e}$ is always within $\kappa''$ of $\mathbf{h}^{\infty,1}$. In the previous statement and below, the constant $\kappa'' \in \R_{>0}$ is universal:
\begin{align}
\mathbf{z}_{T,X}^{\infty,1} \mathrm{e}^{-\kappa''} \ \leq \ \mathbf{z}_{T,X}^{\infty,\e} \ \leq \ \mathbf{z}_{T,X}^{\infty,\e} \ \leq \ \mathbf{z}_{T,X}^{\infty,1} \mathrm{e}^{\kappa''}. \label{eq:CompactSHE3}
\end{align}
To propagate \eqref{eq:CompactSHE3} globally in time, by linearity of the SHE, the space-time fields $\mathbf{z}^{\infty,1}\mathrm{e}^{\pm\kappa''}$ are themselves solutions to the SHE with appropriate initial data setting $T = 0$. In particular, the comparison principle of Theorem 1 in \cite{Mu} allows us to propagate \eqref{eq:CompactSHE3} for all $T \in \R_{\geq 0}$, for all $X \in \mathbb{T}$, and for all $\e \in [0,1]$ simultaneously with probability 1. Thus, we have
\begin{align}
\mathscr{E}_{\e,\delta,\Phi,\mathscr{C}}^{1} \ &\supseteq \ \left\{ \mathrm{e}^{-\kappa'\mathscr{C}} \ \leq \ \inf_{S\in[0,T]}\inf_{X\in\mathbb{T}} \mathbf{z}_{S,X}^{\infty,1} \mathrm{e}^{-\kappa''} \ \leq \ \sup_{S\in[0,T]}\sup_{X\in\mathbb{T}} \mathbf{z}_{S,X}^{\infty,1} \mathrm{e}^{\kappa''} \ \leq \ \mathrm{e}^{\kappa'\mathscr{C}} \right\},
\end{align}
so translating this in terms of probabilities of complements, we get via union bound the following estimate in which $\kappa'''=\kappa'+\kappa''$:
\begin{align}
\mathbf{P}\left(\left(\mathscr{E}_{\e,\delta,\Phi,\mathscr{C}}^{1}\right)^{C}\right) \ &\leq \ \mathbf{P}\left(\inf_{S\in[0,T]}\inf_{X\in\mathbb{T}} \mathbf{z}_{S,X}^{\infty,1} \mathrm{e}^{-\kappa'''} \ \leq \ \mathrm{e}^{-\mathscr{C}}\right) \ + \ \mathbf{P}\left(\sup_{S\in[0,T]}\sup_{X\in\mathbb{T}} \mathbf{z}_{S,X}^{\infty,1} \mathrm{e}^{\kappa'''} \ \geq \ \mathrm{e}^{\mathscr{C}}\right). \label{eq:CompactSHE4}
\end{align}
We observe first solutions to the SHE are jointly continuous in space-time with probability 1, as long as the given initial data is continuous with probability 1, for instance. This follows by standard Holder regularity estimates via the Kolmogorov criterion; see \cite{BG} as a reference for such regularity estimate. To use this feature of SHE, first observe $\mathbf{h}_{0,\bullet}^{\infty,1}$ is continuous with probability 1, because it is the Brownian motion conditioned on finitely many events of uniformly positive probability by construction, of non-degenerate $\sigma$-dependent diffusion coefficient. Thus, both of $\mathbf{z}_{0,\bullet}^{\infty,1}\mathrm{e}^{\pm\kappa'''}$ are continuous with probability 1, which therefore implies $\mathbf{z}_{\bullet,\bullet}^{\infty,1}\mathrm{e}^{\pm\kappa'''}$ are jointly Holder continuous in space-time with probability 1, so that provided any $\delta' \in \R_{>0}$, we may choose $\mathscr{C} \in \R_{>0}$ sufficiently large depending only on $\delta' \in \R_{>0}$ such that
\begin{align}
\mathbf{P}\left(\inf_{S\in[0,T]}\inf_{X\in\mathbb{T}} \mathbf{z}_{S,X}^{\infty,1} \mathrm{e}^{-\kappa'''} \ \leq \ \mathrm{e}^{-\mathscr{C}}\right) \ + \ \mathbf{P}\left(\sup_{S\in[0,T]}\sup_{X\in\mathbb{T}} \mathbf{z}_{S,X}^{\infty,1} \mathrm{e}^{\kappa'''} \ \geq \ \mathrm{e}^{\mathscr{C}}\right) \ &\leq \ \delta'. \label{eq:CompactSHE5}
\end{align}
We now estimate the second probability on the RHS of \eqref{eq:CompactSHE2} via space-time pointwise moments of the following difference:
\begin{align}
\mathbf{x}_{T,X}^{\infty,\e} \ \overset{\bullet}= \ \mathbf{z}_{T,X}^{\infty} - \mathbf{z}_{T,X}^{\infty,\e}.
\end{align}
By linearity of the SHE $\mathbf{x}_{\bullet,\bullet}^{\infty,\e}$ solves an identical SHE as well. In particular, courtesy of the Ito isometry for the space-time white noise as well as convexity of the heat propagator $\mathrm{e}^{2^{-1}\alpha T\Delta}$ with heat kernel $\mathscr{U}_{T,X}$, we have the following in which $\mathscr{L}^{p}_{\omega}$ denotes the $p$-norm with respect to the underlying randomness of the coupled initial conditions and the randomness of the noise:
\begin{align}
\|\mathbf{x}_{T,X}^{\infty,\e}\|_{\mathscr{L}^{2}_{\omega}}^{2} \ &\lesssim \ \int_{\mathbb{T}} \mathscr{U}_{T,X-Y} \|\mathbf{x}_{0,Y}^{\infty,\e}\|_{\mathscr{L}^{2}_{\omega}}^{2} \ \d Y \ + \ \int_{0}^{T} \int_{\mathbb{T}} \left|\mathscr{U}_{T-S,X-Y}\right|^{2} \|\mathbf{x}_{S,Y}^{\infty,\e}\|_{\mathscr{L}^{2}_{\omega}}^{2} \ \d Y \d S \\
&\lesssim \ \int_{\mathbb{T}} \mathscr{U}_{T,X-Y} \|\mathbf{x}_{0,Y}^{\infty,\e}\|_{\mathscr{L}^{2}_{\omega}}^{2} \ \d Y \ + \ \int_{0}^{T} \int_{\mathbb{T}} |T-S|^{-\frac12} \mathscr{U}_{T-S,X-Y} \|\mathbf{x}_{S,Y}^{\infty,\e}\|_{\mathscr{L}^{2}_{\omega}}^{2} \ \d Y \d S \\
&\leq \ \sup_{X \in \mathbb{T}} \|\mathbf{x}_{0,Y}^{\infty,\e}\|_{\mathscr{L}^{2}_{\omega}}^{2} \ + \ \int_{0}^{T} |T-S|^{-\frac12} \sup_{X \in \mathbb{T}} \|\mathbf{x}_{S,Y}^{\infty,\e}\|_{\mathscr{L}^{2}_{\omega}}^{2} \ \d S.
\end{align}
Above, we have used the heat kernel estimate $\mathscr{U}_{T,X} \lesssim T^{-\frac12}$ for order 1 times with universal implied constant, which is classical. As this final upper bound is uniform in the spatial coordinate, another application of the singular Gronwall inequality gives the following upper bound with implied constant depending only on $T \in \R_{\geq 0}$:
\begin{align}
\sup_{X\in\mathbb{T}}\|\mathbf{x}_{T,X}^{\infty,\e}\|_{\mathscr{L}^{2}_{\omega}}^{2} \ &\lesssim_{T} \ \sup_{X\in\mathbb{T}}\|\mathbf{x}_{0,X}^{\infty,\e}\|_{\mathscr{L}^{2}_{\omega}}^{2}
\end{align}
By the height function estimate from Corollary \ref{corollary:ASA} after taking the large-$N$ limit and by the deterministic uniformly bounded nature of $\mathbf{h}^{\infty}$, elementary calculations show the RHS vanishes as $\e \to 0$; indeed, comparison of $\mathbf{h}^{\infty}=\mathbf{h}^{\infty,0}$ and $\mathbf{h}^{\infty,\e}$ from Corollary \ref{corollary:ASA} extends to that of their exponentials because they are all uniformly bounded above and below since $\mathbf{h}^{\infty}$ is deterministic. Thus, uniformly on compact space-time sets, we get
\begin{align}
\lim_{\e \to 0^{+}} \sup_{X\in\mathbb{T}}\|\mathbf{x}_{T,X}^{\infty,\e}\|_{\mathscr{L}^{2}_{\omega}}^{2} \ = \ 0.
\end{align}
Thus, combining the previous vanishing with the Chebyshev inequality, recalling definition of $\mathbf{x}^{\infty,\e}$ we have
\begin{align}
\limsup_{\e\to0^{+}} \ \mathbf{P}\left(\left(\mathscr{E}_{\e,\delta,\Phi,\mathscr{C}'}^{2}\right)^{C}\right) \ &\leq \ \limsup_{\e \to 0^{+}} \ \|\int_{\mathbb{T}} | \mathbf{z}_{T,X}^{\infty} - \mathbf{z}_{T,X}^{\infty,\e} | \ \d X \|_{\mathscr{L}^{2}_{\omega}}^{2} \cdot \delta^{-2} |\mathscr{C}'|^{2} \|\Phi\|_{\mathscr{L}^{\infty}_{X}}^{2} \\
&\leq \limsup_{\e \to 0^{+}} \ \int_{\mathbb{T}} \| \mathbf{z}_{T,X}^{\infty} - \mathbf{z}_{T,X}^{\infty,\e}\|_{\mathscr{L}^{2}_{\omega}}^{2} \ \d X \cdot \delta^{-2} |\mathscr{C}'|^{2} \|\Phi\|_{\mathscr{L}^{\infty}_{X}}^{2} \\
&= \ 0. \label{eq:CompactSHE6}
\end{align}
Ultimately, combining \eqref{eq:CompactSHE2}, \eqref{eq:CompactSHE4}, \eqref{eq:CompactSHE5}, and \eqref{eq:CompactSHE6}, we have the following upper bound for any $\delta' \in \R_{>0}$:
\begin{align}
\limsup_{\e\to0^{+}} \ \mathbf{P}\left(\mathscr{E}_{\e,\delta,\Phi}\right) \ &\leq \ \delta'.
\end{align}
As $\delta' \in \R_{>0}$, is arbitrary, we deduce $\lim_{\e\to0^{+}}\mathbf{P}[\mathscr{E}_{\e,\delta,\Phi}] = 0$, thereby completing the proof.
\end{proof}
%
%
%
\section{Proofs of Theorem \ref{theorem:SBEZRP} and Theorem \ref{theorem:ASHEZRP}}
\subsection{Proof of Theorem \ref{theorem:SBEZRP}}
The proof of Theorem \ref{theorem:SBEZRP} consists of tightness and identifying subsequential limit points. Both depend on the following summation-by-parts identity. It will be convenient to establish the following notation. First, let us assume that the speed $c_{N,\mathfrak{a},\gamma,\sigma}$ from Definition \ref{definition:intro3} is equal to zero; this is certainly not necessarily true, but the exact value of $c_{N,\mathfrak{a},\gamma,\sigma}$ will not be important for the rest of our analysis and makes our notation and presentation more convenient.
\begin{notation}
Provided any $\e \in \R_{\geq 0}$ and smooth test function $\Phi \in \mathscr{C}^{\infty}(\mathbb{T})$, we define the process
\begin{align}
T \ \rightsquigarrow \ \mathbf{Y}_{T}^{N,\e}(\Phi) \ \overset{\bullet}= \ N^{-\frac12} {\sum}_{X \in \mathbb{T}_{N}} \left( \eta_{T,X}^{N,\e} - \sigma \right) \cdot \Phi(N^{-1} X),
\end{align}
where $\eta_{\bullet}^{N,\e}$ denotes the particle configuration obtained by sampling the initial condition via the probability measure $\mathbf{P}_{N,\e}$ from \emph{Corollary \ref{corollary:ASA}} jointly independently over all $\e \in \R_{\geq 0}$. Moreover, we couple all of the initial conditions for the particle dynamics via the basic coupling; see \emph{Remark \ref{remark:InfCoupling}}. Under this basic coupling, define the following given $\e \in \R_{>0}$ and $\Phi \in \mathscr{C}^{\infty}(\mathbb{T})$:
\begin{align}
T \ \rightsquigarrow \ \bar{\mathbf{Y}}_{T}^{N,\e}(\Phi) \ \overset{\bullet}= \ \mathbf{Y}_{T}^{N,0}(\Phi) - \mathbf{Y}_{T}^{N,\e}(\Phi).
\end{align}
\end{notation}
\begin{lemma}\label{lemma:SBP}
Consider any $\e \in \R_{\geq 0}$ and any smooth test function $\Phi \in \mathscr{C}^{\infty}(\mathbb{T})$. We have 
\begin{align}
\mathbf{Y}_{T}^{N,\e}(\Phi) \ = \ N^{-1}{\sum}_{X \in \mathbb{T}_{N}} \grad_{-1}^{N} \Phi(N^{-1} X) \cdot \mathbf{h}_{T,X}^{N,\e} \ + \ \grad_{N-1}\mathbf{h}_{T,0}^{N,\e} \cdot \Phi(1).
\end{align}
\end{lemma}
\begin{proof}
By definition of $\mathbf{h}^{N,\e}$ and $\mathbf{Y}^{N,\e}$, we have the following which resembles $\mathbf{Y}^{N,\e} \ ``=" \ \grad\mathbf{h}^{N,\e}$. We just have to take extra care at the endpoint $N\in\mathbb{T}_{N}$ and note that $\grad_{+}^{N}\mathbf{h}^{N,\e}_{T,N}=N^{1/2}\sum_{X\in\mathbb{T}_{N}}(\eta_{T,X}^{N,\e}-\sigma)$; see Definition \ref{definition:HF}. This is actually where the second term in the proposed identity in Lemma \ref{lemma:SBP} comes from.
\begin{align*}
\mathbf{Y}_{T}^{N,\e}(\Phi) \ &= \ N^{-1} {\sum}_{X \in \mathbb{T}_{N} \setminus \{N\}} \grad_{+1}^{N} \mathbf{h}_{T,X}^{N,\e} \cdot \Phi(N^{-1} X) \ + \ N^{-\frac12} \left(\eta_{T,N}^{N,\e} - \sigma\right) \cdot \Phi(1) \\
&= \ N^{-1} {\sum}_{X \in \mathbb{T}_{N}} \grad_{+1}^{N} \mathbf{h}_{T,X}^{N,\e} \cdot \Phi(N^{-1} X) \ - \ N^{-\frac12} {\sum}_{X \in \mathbb{T}_{N}} \left( \eta_{T,X}^{N,\e} - \sigma \right) \cdot \Phi(1) \ + \ N^{-\frac12} \left( \eta_{T,N}^{N,\e} - \sigma \right) \cdot \Phi(1) \\
&= \ N^{-1} {\sum}_{X \in \mathbb{T}_{N}} \grad_{+1}^{N} \mathbf{h}_{T,X}^{N,\e} \cdot \Phi(N^{-1} X) \ - \ N^{-\frac12} {\sum}_{X \in \mathbb{T}_{N} \setminus \{N\}} \left( \eta_{T,X}^{N,\e} - \sigma \right) \cdot \Phi(1).
\end{align*}
The claimed identity follows from summation-by-parts on the torus $\mathbb{T}_{N}$ for the first quantity on the RHS of the last display. Indeed, the second term in the last line is the discrete gradient $\grad_{N-1}\mathbf{h}_{T,0}^{N,\e}$ by definition of the height function. This completes the proof.
\end{proof}
As consequence of this summation-by-parts-type identity, we obtain the following comparison between $\mathbf{Y}^{N,\e}$ and $\mathbf{Y}^{N,0}$ for all $\e \in \R_{>0}$. Actually, the utility behind Lemma \ref{lemma:SBP} for proving Theorem \ref{theorem:SBEZRP} is exactly through providing this upcoming result, which will both be important towards establishing tightness of the fluctuation field $\mathbf{Y}^{N,0}$ and identifying limit points as the distributional derivative of the Cole-Hopf solution to the KPZ equation.
\begin{corollary}\label{corollary:ComparisonFF}
Consider any $\e \in \R_{>0}$ and any smooth test function $\Phi \in \mathscr{C}^{\infty}(\mathbb{T})$. Deterministically in $T \in \R_{\geq 0}$, we have
\begin{align}
\bar{\mathbf{Y}}_{T}^{N,\e}(\Phi) \ = \ N^{-1} {\sum}_{X\in\mathbb{T}_{N}} \grad_{-1}^{N} \Phi(N^{-1}X) \cdot \left(\mathbf{h}_{T,X}^{N,0} - \mathbf{h}_{T,X}^{N,\e}\right) \ + \ \left(\grad_{N-1}\mathbf{h}_{T,0}^{N,0} - \grad_{N-1}\mathbf{h}_{T,0}^{N,\e} \right) \cdot \Phi(1). \label{eq:ComparisonFF00}
\end{align}
Thus, we have the following estimate uniformly in $T \in \R_{\geq 0}$ with probability at least $1 - \delta_{N}$, where $\delta_{N} \to_{N \to \infty} 0$:
\begin{align}
\|\bar{\mathbf{Y}}_{T}^{N,\e}(\Phi)\|_{\mathscr{L}^{\infty}_{T}} \ &\lesssim_{\Phi} \ \e. \label{eq:ComparisonFF01}
\end{align}
To be completely explicit, the implied constant is only allowed to depend on a finite set of norms of $\Phi \in \mathscr{C}^{\infty}(\mathbb{T})$.
\end{corollary}
\begin{proof}
The identity \eqref{eq:ComparisonFF00} is straightforward to establish given Lemma \ref{lemma:SBP}, so it remains to get \eqref{eq:ComparisonFF01}. We first note the following consequence of a number of things including the fundamental theorem of calculus, which holds uniformly in $X \in \mathbb{T}_{N}$:
\begin{align}
|\grad_{-1}^{N}\Phi(N^{-1}X)| \ &\leq \ \| \partial_{X} \Phi\|_{\mathscr{L}^{\infty}_{X}}.
\end{align}
Before we proceed, we clarify that all statements to follow in this argument that we say hold with probability at least $1 - \delta_{N}$ hold simultaneously on the same event with probability at least $1 - \delta_{N}$. We will only make a uniformly bounded number of such statements, so the final conclusion holds with probability $1 - \mathscr{O}(\delta_{N})$ by union bound for the events on which our statements fail.

Courtesy of the basic coupling and Lemma \ref{lemma:ComparisonMicro} and the above derivative bound, with probability at least $1 - \delta_{N}$ we get
\begin{align}
\| N^{-1} {\sum}_{X\in\mathbb{T}_{N}} \grad_{-1}^{N} \Phi(N^{-1}X) \cdot \left(\mathbf{h}_{T,X}^{N,0} - \mathbf{h}_{T,X}^{N,\e}\right) \|_{\mathscr{L}^{\infty}_{T}} \ &\leq \ \|\partial_{X} \Phi\|_{\mathscr{L}^{\infty}_{X}} \|\mathbf{h}_{T,X}^{N,0}-\mathbf{h}_{T,X}^{N,\e}\|_{\mathscr{L}^{\infty}_{T}\mathscr{L}^{\infty}_{X}} \ \lesssim \ \|\partial_{X} \Phi\|_{\mathscr{L}^{\infty}_{X}} \e. \label{eq:ComparisonFF1}
\end{align}
The implied constant above is universal. It remains to address the gradient difference of $\mathbf{h}^{N,0}$ and $\mathbf{h}^{N,\e}$ on the RHS of \eqref{eq:ComparisonFF00}. For this, we employ the $\e$-closeness of the height functions $\mathbf{h}^{N,0}$ and $\mathbf{h}^{N,\e}$ in Lemma \ref{lemma:ComparisonMicro}. We first write
\begin{align}
|\grad_{N-1}\mathbf{h}_{T,0}^{N,0} - \grad_{N-1}\mathbf{h}_{T,0}^{N,\e}| \ &\leq \ |\mathbf{h}_{T,N-1}^{N,0}-\mathbf{h}_{T,N-1}^{N,\e}| + |\mathbf{h}_{T,0}^{N,0}-\mathbf{h}_{T,0}^{N,\e}|.
\end{align}
By Lemma \ref{lemma:ComparisonMicro}, from the previous estimate we get with a universal implied constant
\begin{align}
|\grad_{N-1}\mathbf{h}_{T,0}^{N,0} - \grad_{N-1}\mathbf{h}_{T,0}^{N,\e}| \cdot |\Phi(1)| \ &\lesssim_{\Phi} \ \e.
\end{align}
Combining this previous estimate with \eqref{eq:ComparisonFF00} and \eqref{eq:ComparisonFF1} completes the proof.
\end{proof}
\begin{proof}[Proof of \emph{Theorem \ref{theorem:SBEZRP}}]
As noted earlier, we prove tightness and identify limit points. Take any finite time-horizon $T \in \R_{\geq 0}$.
\begin{itemize}[leftmargin=*]
\item We first show tightness of $\mathbf{Y}^{N,0}: [0,T] \to \mathscr{D}(\mathbb{T})$. By the Mitoma criterion in \cite{M83}, it suffices to show tightness of the $\R$-valued process $\mathbf{Y}^{N,0}(\Phi): [0,T] \to \R$ for any fixed $\Phi \in \mathscr{C}^{\infty}(\mathbb{T})$. By convergence in Corollary \ref{corollary:ASA}, the process $\mathbf{Y}^{N,\e}(\Phi): [0,T] \to \R$ is tight in the large-$N$ limit for any $\e \in \R_{>0}$ fixed. As $\e \in \R_{>0}$ may be taken arbitrarily small, combining this with Corollary \ref{corollary:ComparisonFF} and Lemma \ref{lemma:ApproxTight}, with choices $\mathfrak{C}^{N,\e}=\mathbf{Y}^{N,\e}(\Phi)$ and $\mathfrak{A}^{N}=\mathbf{Y}^{N,0}(\Phi)$, yields tightness of $\mathbf{Y}^{N,0}(\Phi)$.
\item We proceed to identify subsequential limit points. Consider any $\Phi \in \mathscr{C}^{\infty}(\mathbb{T})$ and any limit point $\mathbf{Y}^{\infty}$. Given arbitrarily small but $N$-independent $\e \in \R_{>0}$, Corollary \ref{corollary:ComparisonFF} provides the following estimate for any limit point $\mathbf{Y}^{\infty,\e} = \lim_{N\to\infty}\mathbf{Y}^{N,\e}$ and $T \in \R_{\geq 0}$. Note in the following the implied constant is independent of $T\in\R_{\geq0}$:
\begin{align}
|\mathbf{Y}_{T}^{\infty}(\Phi) - \mathbf{Y}_{T}^{\infty,\e}(\Phi) | \ \lesssim_{\Phi} \ \e.
\end{align}
In particular, we get $\mathbf{Y}_{T}^{\infty}(\Phi) = \lim_{\e\to0^{+}}\mathbf{Y}_{T}^{\infty,\e}(\Phi)$. However, applying Theorem 2.13 from \cite{GP}, we also have the distributional identity $\mathbf{Y}_{T}^{\infty,\e}(\Phi) = -\beta\alpha^{-1} \langle\log\mathbf{z}_{T,X}^{\infty,\e},\partial_{X}\Phi\rangle$, where the additional negative sign follows by integrating by parts and moving the $\partial_{X}$ derivative onto $\Phi$. Therefore, by Lemma \ref{lemma:CompactSHE}, we get the following convergence in probability:
\begin{align}
\mathbf{Y}_{T}^{\infty}(\Phi) \ &= \ \lim_{\e\to0^{+}} -\beta\alpha^{-1} \int_{\mathbb{T}} \partial_{X} \Phi(X) \cdot \log \mathbf{z}_{T,X}^{\infty,\e} \ \d X \ = \ -\beta\alpha^{-1} \int_{\mathbb{T}} \partial_{X} \Phi(X) \cdot \log \mathbf{z}_{T,X}^{\infty} \ \d X.
\end{align}
Thus, pointwise-in-time, with probability 1 we get the distributional identity $\mathbf{Y}^{\infty} = \beta\alpha^{-1}\partial_{X} \log \mathbf{z}^{\infty}$. As both distributions are continuous with respect to the time-variable, we may extend this identity globally in time with probability 1. The continuity with respect to the time-variable follows from vanishing in the large-$N$ limit of jumps in $\mathbf{Y}^{N,0}$ and $\mathbf{Y}^{N,\e}$ and the fact that any limit point in the Skorokhod space $\mathscr{D}(\R_{\geq0},\mathscr{D}(\mathbb{T}))$ whose jumps are all of size zero must be a continuous function; see \cite{Bil} for details about the Skorokhod space and its relationship with the space $\mathscr{C}(\R_{\geq0},\mathscr{D}(\mathbb{T}))$ of continuous paths.
\end{itemize}
This completes the proof.
\end{proof}
\begin{proof}[Proof of \emph{Theorem \ref{theorem:ASHEZRP}}]
The proof of Theorem \ref{theorem:ASHEZRP} is identical to the proof of Theorem \ref{theorem:SBEZRP} but instead of the SPDE result Lemma \ref{lemma:CompactSHE}, we apply the following elementary convergence result for $\mathrm{EW}$.
\begin{lemma}\label{lemma:CompactASHE}
Retain the setting of \emph{Lemma \ref{lemma:CompactSHE}}, except all ``Cole-Hopf solutions" to the \emph{KPZ equation} are replaced by solutions to the $\mathrm{EW}$. For any smooth test function $\Phi \in \mathscr{C}^{\infty}(\mathbb{T})$ along with any time $T \in \R_{\geq 0}$, we have the convergence in probability
\begin{align}
\int_{\mathbb{T}} \Phi(X) \cdot \mathbf{h}_{T,X}^{\infty,\e} \ \d X \ \longrightarrow_{\e \to 0^{+}} \ \int_{\mathbb{T}} \Phi(X) \cdot \mathbf{h}_{T,X}^{\infty} \ \d X.
\end{align}
The solutions $\mathbf{h}^{\infty}$ and $\{\mathbf{h}^{\infty,\e}\}_{\e\in\R_{>0}}$ are coupled by independently sampling initial data and coupling to the same noise.
\end{lemma}
This completes the proof of Theorem \ref{theorem:ASHEZRP}.
\end{proof}
\appendix
\section{Technical Lemmas}
The first technical result we prove in this appendix section is a result from stochastic calculus known as a support theorem. It is the feature of Wiener measure that gives open balls positive probability.
\begin{lemma}\label{lemma:SupportTheorem}
Suppose $\mathfrak{f} \in \mathscr{C}([0,1],\R)$ is a deterministic continuous function satisfying $\mathfrak{f}(0) = \mathfrak{f}_{-}$ and $\mathfrak{f}(1) = \mathfrak{f}_{+}$. Moreover, suppose $\{\mathfrak{B}(x)\}_{x\in[0,1]}$ is a Brownian motion satisfying $\mathfrak{B}(0) = \mathfrak{f}_{-}$ conditioned on the event
\begin{align}
\mathscr{E}_{1,\e} \ \overset{\bullet}= \ \left\{ |\mathfrak{B}(1)-\mathfrak{f}_{+}| \leq \e \right\}.
\end{align}
Given any $\e' \in \R_{>0}$ satisfying $\e' \leq \e$, consider the event $\mathscr{E}_{\e'} \overset{\bullet}=\{\sup_{x\in[0,1]}|\mathfrak{B}(x) - \mathfrak{f}(x)|\leq\e'\}$. For any such $\e' \in \R_{>0}$, we have
\begin{align}
\mathbf{P}\left(\mathscr{E}_{\e'}\right) \ \gtrsim_{\e',\e,\mathfrak{f}} \ 1.
\end{align}
\end{lemma}
\begin{proof}
As $\mathbf{P}[\mathscr{E}_{1,\e}] \gtrsim_{\e} 1$ and $\mathscr{E}_{\e'} \subseteq \mathscr{E}_{1,\e}$, it suffices to assume $\mathfrak{B}$ is instead the Brownian motion free of any conditioning. This can be checked with definition of conditional probability given $\mathbf{P}[\mathscr{E}_{1,\e}] \gtrsim_{\e} 1$, for example. Moreover, we may assume $\mathfrak{f} \in \mathscr{C}([0,1],\R)$ is actually smooth. Indeed, if $\wt{\mathfrak{f}}$ is a smooth function satisfying $\sup_{x\in[0,1]}|\mathfrak{f}(x)-\wt{\mathfrak{f}}(x)| \leq \e'$, for example a smooth convolution of $\e',\mathfrak{f}$-dependent scale, then we have the following triangle inequality
\begin{align}
{\sup}_{x\in[0,1]}|\mathfrak{B}(x) - \mathfrak{f}(x)| \ &\leq \ {\sup}_{x\in[0,1]}|\mathfrak{B}(x) - \wt{\mathfrak{f}}(x)| \ + \ {\sup}_{x\in[0,1]}|\wt{\mathfrak{f}}(x) - \mathfrak{f}(x)|
\end{align}
from which the result for $\wt{\mathfrak{f}}$ with parameter $\e' \in \R_{>0}$ would give the result for $\mathfrak{f}$ with parameter $2\e' \in \R_{>0}$. However, in this reduction, we emphasize that the regularity of the smooth function $\wt{\mathfrak{f}}$ now depends on $\e,\e' \in \R_{>0}$ and the modulus of continuity of $\mathfrak{f}$. Because $\mathfrak{B}$ is a standard Brownian motion and $\mathfrak{f}$ is smooth, courtesy of the Girsanov theorem and the Novikov condition, the law of the standard Brownian motion is absolutely continuous with respect to the law of the Ito diffusion $\mathfrak{B} - \mathfrak{f}$, and the corresponding Radon-Nikodym derivative is an $\mathscr{L}^{2}$-martingale with respect to the path-space probability measure of the process $\mathfrak{B}-\mathfrak{f}$ whose $\mathscr{L}^{2}$-norm depends only on analytic data of $\mathfrak{f}$, which we recall depends on $\e,\e' \in \R_{>0}$. Ultimately, by the Cauchy-Schwarz inequality for functionals/events on the path space $\mathscr{C}([0,1],\R)$, we have
\begin{align}
\mathbf{P}({\sup}_{x\in[0,1]}|\mathfrak{B}(x)|\leq \e')^{2} \ &\lesssim_{\mathfrak{f},\e,\e'} \ \mathbf{P}\left(\mathscr{E}_{\e'}\right).
\end{align}
Upon bounding the probability on the LHS below via the reflection principle, this completes the proof.
\end{proof}
The second result we establish concerns the topology of the Skorokhod space $\mathscr{D}(\R_{\geq 0},\R)$. In short, the content behind the following result is that uniform approximations of tight sequences are also tight.
\begin{lemma}\label{lemma:ApproxTight}
Suppose $\{\mathfrak{A}_{\bullet}^{N}\}_{N \in \Z_{>0}}$ is any family of stochastic processes with its sample paths in $\mathscr{D}(\R_{\geq 0},\R)$. Moreover, suppose we further have a one-parameter family of families $\{\mathfrak{C}_{\bullet}^{N,\e}\}_{N\in\Z_{>0}}$ of stochastic processes indexed by $\e \in \R_{>0}$ with sample paths in $\mathscr{D}(\R_{\geq 0},\R)$. Lastly, suppose the following conditions hold:
\begin{itemize}[leftmargin=*]
\item For each fixed $\e \in \R_{>0}$, the family $\{\mathfrak{C}_{\bullet}^{N,\e}\}_{N\in\Z_{>0}}$ is tight with respect to the Skorokhod topology on $\mathscr{D}(\R_{\geq 0},\R)$.
\item For each fixed $\e \in \R_{>0}$, we have the following estimate globally in time in which $\delta_{N}\to_{N\to\infty}0$ given any fixed positive $\e$:
\begin{align}
\mathbf{P}\left({\sup}_{N \in \Z_{>0}} \| \mathfrak{A}_{T}^{N} - \mathfrak{C}_{T}^{N,\e} \|_{\mathscr{L}^{\infty}_{T\geq 0}} \lesssim \e\right) \ \geq \ 1 - \delta_{N}.\label{eq:ApproxTight0}
\end{align}
\end{itemize}
Then, the family $\{\mathfrak{A}_{\bullet}^{N}\}_{N\in\Z_{>0}}$ is itself tight with respect to the Skorokhod topology on $\mathscr{D}(\R_{\geq 0},\R)$.
\end{lemma}
\begin{proof}
First, note that it suffices to replace the non-compact time-interval by $[0,T_{f}] \subseteq \R_{\geq 0}$ where $T_{f} \in \R_{\geq 0}$ is some generic time-horizon. In this spirit, for this proof let us define $\|\|_{\mathscr{L}^{\infty}_{T_{f}}}=\sup_{T\in[0,T_{f}]}$. Second, observe the following probability estimate courtesy of the triangle inequality for $\|\|_{\mathscr{L}^{\infty}_{T_{f}}}$, in which $\mathscr{C} \in \R_{>0}$ is a large constant:
\begin{align}
\mathbf{P}\left(\| \mathfrak{A}_{T}^{N} \|_{\mathscr{L}^{\infty}_{T_{f}}} \geq \mathscr{C} \right) \ &\leq \ \mathbf{P}\left(\| \mathfrak{C}_{T}^{N,\e} \|_{\mathscr{L}^{\infty}_{T_{f}}} + \|\mathfrak{A}_{T}^{N} - \mathfrak{C}_{T}^{N,\e} \|_{\mathscr{L}^{\infty}_{T_{f}}} \geq \mathscr{C} \right).
\end{align}
Recalling the assumption \eqref{eq:ApproxTight0}, from the previous upper bound we deduce, for $\e \in \R_{>0}$ sufficiently small,
\begin{align}
\mathbf{P}\left(\| \mathfrak{A}_{T}^{N} \|_{\mathscr{L}^{\infty}_{T_{f}}} \geq \mathscr{C} \right) \ &\leq \ \mathbf{P}\left(\| \mathfrak{C}_{T}^{N,\e} \|_{\mathscr{L}^{\infty}_{T_{f}}} \geq \mathscr{C} - \e \right) \ + \ \delta_{N},
\end{align}
from which we have, if $\e \in \R_{>0}$ is fixed,
\begin{align}
\limsup_{\mathscr{C} \to \infty}\limsup_{N \in \Z_{>0}} \mathbf{P}\left(\| \mathfrak{A}_{T}^{N} \|_{\mathscr{L}^{\infty}_{T_{f}}} \geq \mathscr{C} \right) \ &\leq \ \limsup_{\mathscr{C} \to \infty}\limsup_{N \in \Z_{>0}} \mathbf{P}\left(\| \mathfrak{C}_{T}^{N,\e} \|_{\mathscr{L}^{\infty}_{T_{f}}} \geq \mathscr{C} \right) \ = \ 0. \label{eq:ApproxTight1}
\end{align}
The last bound \eqref{eq:ApproxTight1} follows from the assumed tightness of $\{\mathfrak{C}^{N,\e}_{\bullet}\}_{N\in\Z_{>0}}$ with $\e \in \R_{>0}$ fixed and Theorem 16.8 in \cite{Bil}.

Before we proceed, we use the family of regularity metrics $\mathscr{W}(\bullet,\delta): \mathscr{D}([0,T_{f}],\R) \to \R$ parameterized by $\delta \in \R_{\geq0}$ from Section 16 in \cite{Bil} without defining it explicitly in this article. All we will need is that Theorem 16.8 of \cite{Bil} tells us that given \eqref{eq:ApproxTight1}, tightness of $\mathfrak{A}^{N}$ will follow from control on this regularity metric for $\mathfrak{A}^{N}$. With probability at least $1 - \delta_{N}$ we first get
\begin{align}
\mathscr{W}(\mathfrak{A}^{N},\delta) \ &\leq \ \mathscr{W}(\mathfrak{C}^{N,\e},\delta) \ + \ \mathscr{O}(\e),
\end{align}
where the implied constant in the big-Oh constant is universal outside dependence in the universal implied constant in \eqref{eq:ApproxTight0}. Indeed, this regularity metric estimate is a direct consequence of the definition of $\mathscr{W}(\bullet,\delta)$, namely its triangle inequality as a metric, along with \eqref{eq:ApproxTight0}. Thus, given any $\e' \in \R_{>0}$, let us choose $\e \lesssim_{\e'} 1$ sufficiently small depending only on $\e' \in \R_{>0}$. We then have
\begin{align}
\limsup_{\delta \to 0} \limsup_{N \to \infty} \mathbf{P}\left(\mathscr{W}(\mathfrak{A}^{N},\delta) \geq \e' \right) \ &\leq \ \limsup_{\delta \to 0} \limsup_{N \to \infty} \mathbf{P}\left(\mathscr{W}(\mathfrak{C}^{N,\e},\delta) + \mathscr{O}(\e) \geq \e' \right) \\
&\leq \ \limsup_{\delta \to 0} \limsup_{N \to \infty} \mathbf{P}\left(\mathscr{W}(\mathfrak{C}^{N,\e},\delta) \geq  \frac12\e' \right) \\
&=\ 0, \label{eq:ApproxTight2}
\end{align}
where the last identity \eqref{eq:ApproxTight2} follows from Theorem 16.8 in \cite{Bil} and tightness of $\mathfrak{C}^{N,\e}$, again. Combining \eqref{eq:ApproxTight1} and \eqref{eq:ApproxTight2} with Theorem 16.8 in \cite{Bil} now completes the proof.
\end{proof}
\section{Relative Entropy}
We proceed to introduce the following relative entropy functional which serves as a measurement of stability with respect to the grand-canonical ensemble measures $\mathbf{P}_{\alpha}$ that were constructed in Definition \ref{definition:intro1}.
\begin{definition}\label{definition:RE}
For $\mu_{1},\mu_{2}$ probability measures on $\Omega_{\mathbb{T}_{N}}$ such that $\mu_{1} \ll \mu_{2}$, the \emph{relative entropy} of $\mu_{1}$ with respect to $\mu_{2}$ is
\begin{align}
\mathsf{H}_{\mu_{2}}(\mu_{1}) \ \overset{\bullet}= \ \int_{\Omega_{\mathbb{T}_{N}}} \frac{\d\mu_{1}}{\d\mu_{2}} \log \frac{\d\mu_{1}}{\d\mu_{2}} \ \d\mu_{2}. \label{eq:RE} 
\end{align}
\end{definition}


\end{document}